\documentclass[10pt,reqno]{amsart}
\usepackage[scale=0.75, centering, headheight=14pt]{geometry}
\usepackage{fancyhdr}
\usepackage{amssymb,amsmath,enumerate}
\usepackage{graphics}
\usepackage{graphicx}
\usepackage{caption}
\usepackage{float}
\usepackage{subcaption}
\usepackage[linesnumbered,boxruled,commentsnumbered]{algorithm2e}
\usepackage{comment,color}
\usepackage{tikz}
\usetikzlibrary{positioning,arrows.meta}
\usepackage{amsthm}

\usepackage{mathtools}

\usepackage{cleveref}
\usepackage{booktabs}
\theoremstyle{plain}

\theoremstyle{remark}
\newtheorem{remark}{Remark}

\DeclareMathOperator*{\argmin}{arg\,min}

\graphicspath{{./figures/}{.}}

\numberwithin{equation}{section}

\begin{document}


\title[An adaptive Deep Ritz framework for fully nonlinear PDEs]{An adaptive Deep Ritz framework for second-order fully nonlinear partial differential equations}

\author[A. Caboussat]{Alexandre Caboussat}
\address{
Geneva School of Business Administration (HEG-Gen\`eve), 
University of Applied Sciences and Arts Western Switzerland (HES-SO), 
1227 Carouge, Switzerland, 
Email: \texttt{alexandre.caboussat@hesge.ch} 
}

\author[M. Leclercq]{Martin Leclercq}
\address{
Institute of Mathematics, Ecole Polytechnique F\'ed\'erale de Lausanne, 
1015 Lausanne, Switzerland, 
Email: \texttt{martin.leclercq@epfl.ch}
}

\author[A. Peruso]{Anna Peruso}
\address{
Institute of Mathematics, Ecole Polytechnique F\'ed\'erale de Lausanne, 
1015 Lausanne, Switzerland, 
Email: \texttt{anna.peruso@epfl.ch}
and
Geneva School of Business Administration (HEG-Gen\`eve), 
University of Applied Sciences and Arts Western Switzerland (HES-SO), 
1227 Carouge, Switzerland, 
Email: \texttt{anna.peruso@hesge.ch} 
}

\keywords{Fully nonlinear PDEs, Least-squares, Deep Ritz, Monge-Amp\`ere, Optimal Transport}




\begin{abstract}
As an alternative to PINNs, a Deep Ritz framework is proposed to solve fully nonlinear PDEs. A least-squares algorithm is advocated to decouple the nonlinearities from the variational features of several fully nonlinear PDEs. 
A splitting method allows to iteratively solve local nonlinear problems and linear variational problems at each iteration. 
While existing nonlinear solvers are applied to solve for nonlinearities, we propose a novel coupling with a Deep Ritz neural network approach that is well-suited to the variational flavor of the linear variational problems. 
An adaptive sampling strategy for the selection of collocation points is incorporated to increase the efficiency of the algorithm without sacrificing its accuracy. 

Numerical experiments are presented to solve the Dirichlet problem for several fully nonlinear equations, starting with the prototypical Monge-Amp\`ere equation, showing the flexibility of the approach. 
Numerical results are compared with results obtained using a full PINNs approach. 
Finally, numerical experiments are extended to address the optimal transport Monge-Amp\`ere problem with transport boundary conditions. 
\end{abstract}

\maketitle


\section{Introduction}\label{sec:intro}
Fully nonlinear partial differential equations (PDEs) are equations defined by nonlinearities acting on the highest-order derivative of the unknown function. 
They have been studied quite intensively, from both the theoretical perspective, see, e.g., \cite{CaffarelliCabre1995,villani,dephilippis}, and the perspective of numerical methods \cite{FDmonge,feng,brenner,lakkis2013finite}.

We consider here in particular second-order scalar fully nonlinear equations; in this case, the PDE can be generally written as follows: find $u: \Omega \rightarrow \mathbb{R}$ satisfying
\begin{equation*}
F(D^2u, \nabla u, u, x) = 0, \quad x \in \Omega,
\end{equation*}  
where $\Omega \subset \mathbb{R}^d$, and $D^2u(x)$ and $\nabla u(x)$ denote the Hessian and gradient of $u$ at $x \in \Omega$, respectively. 
The operator $F:\mathbb{R}^{d\times d}_{\text{sym}} \times \mathbb{R}^{d} \times \mathbb{R} \times \Omega \to \mathbb{R}$ is assumed to be nonlinear in at least one of the entries of $D^2u$. This equation is typically supplemented with (Dirichlet) boundary conditions. 

Some examples of operators $F$ include the {\em Monge-Ampère equation}, where $F$ is defined as 
\begin{equation}\label{intro_ma}
F(D^2 u, x) = \det(D^2 u) - f(x), 
\end{equation}
and $f(\cdot)$ is a given function. The {\em Pucci’s equation} consists of the example for which $F$ is given by 
\begin{equation}\label{intro_pucci}
F(D^2 u, x) = \alpha \sum \lambda^+(D^2 u) + \sum \lambda^-(D^2 u) - f(x), 
\end{equation}
where $\alpha > 1$ and $\lambda^+(D^2 u)$, $\lambda^-(D^2 u)$ denote the positive and negative eigenvalues of $D^2 u$, respectively and  $f(\cdot)$ is a given function. 
The {\em $\sigma_2$-equation} corresponds to $F$ defined as 
\begin{equation}\label{intro_sigma2}
F(D^2 u, x) = \sigma_2(D^2 u) - f(x), 
\end{equation}
where $\sigma_2(D^2 u)=\sum_{i<j} \lambda_i \lambda_j$ and $\lambda_i$ are the eigenvalues of $D^2 u$ and  $f(\cdot)$ is a given function. 
Finally, the {\em Minkowski equation} consists of $F$ defined as 
\begin{equation}\label{intro_mink}
F(D^2 u, \nabla u, x) =\det(D^2 u) - K(x) (1 + |\nabla u|^2)^{(d+2)/2} , 
\end{equation}
where $K (\cdot)$ is a given positive function. This is a so-called {\em generalized} Monge-Amp\`ere equation.

These equations play a central role in fields as broad as differential geometry, optimal control, mass transportation, geophysical fluid dynamics, meteorology, or general relativity, see, e.g., \cite{dephilippis,feng,trudinger2008monge}. 
Nonetheless, the well-posedness of these problems is subtle: the existence and uniqueness of solutions typically require additional structural assumptions, such as the convexity of the solution or the monotonicity of the operator, see, e.g., \cite{CaffarelliCabre1995}. 
Moreover, their numerical approximation is particularly challenging due to the lack of a weak formulation coming from their non-divergence formulation, which prevents the direct application of standard variational methods commonly used for linear and semilinear PDEs, see, e.g., \cite{feng}.

\medskip

Several numerical strategies have been developed to address these challenges, both including finite differences schemes \cite{FDmonge}, or finite element based methods \cite{feng,brenner,lakkis2013finite}, such as, for instance, based on least-squares methods and variational principles \cite{GlowinskiICIAM2007,caboussat}. 
More recently, with the increasing interest in physics-informed machine learning \cite{RAISSI2019686,karniadakis2021,Kharazmi2019,deepritz,muller,iterative_DR}, deep learning-based methods have emerged as alternatives for approximating solution of fully nonlinear PDEs, see, e.g., \cite{pinnsconvex,NYSTROM,hacking2025neural}. 

\medskip

In this work, we propose a least-squares relaxation framework for second-order fully nonlinear PDEs, inspired by \cite{caboussat}. 
The key idea is to decouple the nonlinearity from the differential operator, and to solve the two problems in an iterative manner. 
Beyond formulating the approach as a unified framework applicable to a broad class of equations and boundary conditions, we employ the Deep Ritz method \cite{deepritz} to solve the variational subproblems arising from the decoupling. 
When dealing with Monge-Amp\`ere equations, for which the convexity of the solution is required, we advocate the underlying use of Input Convex Neural Networks (ICNNs \cite{icnn}) to inherently satisfy the constraint. 
When the convexity of the solution is not required, we use classical (deep) feedforward neural networks. 

\medskip

To further accelerate training, we introduce an adaptive sampling strategy motivated by importance sampling. The local least-squares defect is used as an indicator of where additional collocation points are most useful. Randomly selected seed points define a Voronoi partition of the domain, and the number of samples assigned to each cell is chosen according to the defect evaluated at its seed, following the sampling strategy of \cite{efficient_sampling}. This yields a simple mesh-free adaptive stratification of the collocation points.

\medskip 

We validate the method on several two-dimensional numerical examples for several fully nonlinear PDEs. 
The performance of the algorithm is compared against a pure PINN model, as implemented in \cite{pinnsconvex,NYSTROM} on standard test cases for the Monge-Amp\`ere equation. 
The numerical results show comparable accuracy for simple cases, while the Deep Ritz approach outperforms PINNs for more challenging ones (see \cite{failure_modes} for similar conclusions). 
We further demonstrate the applicability of the approach to the Minkowski \cite{minkowski} and Pucci's problems \cite{pucci}, and extend the treatment of other boundary conditions in order to provide additional results for the optimal transport problem \cite{dephilippis,pinnsconvex,prins}. 
Finally, we show that the proposed adaptive sampling strategy significantly accelerates the learning at the first stages of the splitting algorithm, underlying its importance in the overall design of the algorithm.

\medskip

This article is structured as follows. In \Cref{sec:model}, we describe the PDE problem and decompose it with a least-squares algorithm.  \Cref{sec:ext} details the extensions to several other fully nonlinear equations. 
The Deep Ritz variational algorithm is described in \Cref{sec:deepritz}, together with the underlying neural networks used in the algorithms. 
It also includes the description of the adaptive sampling strategy. 
\Cref{sec:impl} shortly describes the implementation details of the algorithm. 
\Cref{sec:numres,sec:numres_ext} cover a broad range of numerical experiments for Dirichlet problems and optimal transport respectively, before concluding in \Cref{sec:conc}. 


\section{Model background}\label{sec:model}
The framework we consider here has been originally introduced in \cite{caboussat} to approximate smooth solutions to the Dirichlet Monge-Amp\`ere equation. 
Then it has been extended to other fully nonlinear PDEs in various works, see, e.g., \cite{pucci,prins,dimitrios,lsjac}. 
The main idea is to introduce a least-squares formulation of the nonlinear PDE together with an iterative algorithm to decouple the nonlinearity from the highest-order derivatives. The resulting differential problem is to be treated with classical variational methods, while the nonlinearities are solved with mathematical programming techniques. 

More precisely, for a generic operator $F:\mathbb{R}^{d\times d}_{\text{sym}} \to \mathbb{R}$, we first focus on Dirichlet problems of the form: find $u: \Omega \rightarrow \mathbb{R}$ satisfying 
\begin{equation}\label{eq:fnl}
\begin{cases}  
F(D^2u) = 0 \quad &\text{in } \Omega, \\  
u = \phi \quad &\text{on } \partial \Omega,  
\end{cases}  
\end{equation}  
where the explicit dependence on the space variable $x$ has been omitted for readability purposes. 
Then we will discuss extensions to $F=F(D^2u,\nabla u,u)$ (such as, e.g., the Minkowski problem) and to other boundary conditions such as nonlinear Neumann conditions for the optimal transport problem.

Problem (\ref{eq:fnl}) is transformed as follows: let us introduce a matrix field $\mathbf{P}:\Omega\to\mathbb{R}^{d\times d}_{\text{sym}}$. 
With this auxiliary variable the original problem can be rewritten in an equivalent coupled form: find $(u,\mathbf{P})$ satisfying 
\begin{equation}\label{eq:aux}
\begin{cases}
D^2u = \mathbf{P} \quad &\text{in } \Omega, \\  
F(\mathbf{P}) = 0 \quad &\text{in }  \Omega, \\  
u = \phi \quad &\text{on } \partial \Omega. 
\end{cases}
\end{equation}
If $u\in H^2(\Omega)$ is a solution to \eqref{eq:fnl}, then the pair $(u,\mathbf{P})=(u,D^2u)$ satisfies \eqref{eq:aux}. 
We next cast the problem \eqref{eq:aux} as a nonlinearly constrained least-squares minimization problem:
\begin{equation}
\label{eq:leastsq}
(u, \mathbf{P}) = \argmin_{\substack{v\in H^2(\Omega)\cap H^1_\phi(\Omega) \\ \mathbf{Q}\in  L^2(\Omega;\mathbb{R}^{d\times d}_{\text{sym}})}}\left\{J(v,\mathbf{Q}), \quad\text{s.t. } F(\mathbf{Q}) = 0\right\},
\end{equation}  
where 
\begin{equation*}
J(v,\mathbf{Q})= \frac12 \int_\Omega |D^2v - \mathbf{Q}|^2,
\end{equation*}
$|\cdot|$ denotes the Frobenius norm, and $H^1_\phi(\Omega) = \{ u \in H^1 (\Omega) : u = \phi \textrm{ on } \partial \Omega \}$. The consistency with the differential problems \eqref{eq:fnl} and \eqref{eq:aux} is enforced via the nonlinear constraint. 

\medskip

\begin{remark}\label{remark1}
When considering the Dirichlet Monge-Amp\`ere problem, the nonlinearity reads 
$$
F(\mathbf{Q})=\det(\mathbf{Q})-f.
$$
On the other hand, when considering, e.g., the Pucci's equation, the nonlinearity reads
$$
F(\mathbf{Q})= \alpha \sum \lambda^+(\mathbf{Q}) + \sum \lambda^-(\mathbf{Q}) - f,
$$
where $\alpha > 1$ and $\lambda^+(\mathbf{Q})$, $\lambda^-(\mathbf{Q})$ denote the positive and negative eigenvalues of $\mathbf{Q}$.
\end{remark}


\medskip

Since the direct minimization of \eqref{eq:leastsq} is cumbersome, the problem is typically addressed via a \textit{block coordinate descent} approach to decouple the nonlinearity from the differential operator \cite{caboussat,peruso2026}. 
The splitting leads to an iterative scheme in which a constrained minimization problem and a differential problem are solved separately. 
Specifically, given an initial solution $u^0 $, for $n \geq 0$, we seek successively $\mathbf{P}^n $ and $ u^{n+1} $ such that:
\begin{align}
\label{eq:firstmin}
\mathbf{P}^n &= \argmin_{\mathbf{Q}\in  L^2(\Omega; \mathbb{R}^{d\times d}_{\text{sym}})}\left\{J(u^n,\mathbf{Q}), \quad\text{s.t. } F(\mathbf{Q}) = 0\right\},\\
\label{eq:biharmonic}
u^{n+1} &= \argmin_{v\in H^2(\Omega)\cap H^1_\phi(\Omega)}J(v,\mathbf{P}^n),
\end{align}  
In this splitting algorithm the nonlinearity is isolated in the first subproblem \eqref{eq:firstmin}, while the second subproblem \eqref{eq:biharmonic} handles the, linear, variational part of the problem, independently of $F$. However, in its strong formulation, \eqref{eq:biharmonic} is a fourth-order biharmonic equation and typically constitutes the computational bottleneck. 
In \cite{caboussat,perusoestimates} the authors advocate a $\mathbb{P}_1$ finite element scheme and obtain optimal convergence rates for smooth problems on polygonal domains. 

Since \eqref{eq:biharmonic} is given in its variational formulation, in this work we propose using the Deep Ritz method \cite{deepritz} to solve it.

On the other hand, the constrained problem \eqref{eq:firstmin} is nonlinear but local and it can be solved pointwise for any $x\in\Omega$. 
Given $u^n$, for any $x\in \Omega$, we look for the closest matrix to $D^2u^n(x)$ in the Frobenius norm that satisfies the constraint $F(\mathbf{Q}) = 0$. 
Numerically, \eqref{eq:firstmin} is commonly handled by parametrizing $\mathbf{Q}$ (\textit{e.g.} via a factorization or eigenvalue decomposition) and enforcing the constraint through Lagrange multipliers, thereby reformulating it as an unconstrained optimization problem to which standard solvers can be applied. 
Moreover, when the constraint is rotation-invariant (which is the case for many examples), the pointwise optimization reduces to an optimization over the eigenvalues of $\mathbf{Q}$, which significantly reduces the computational cost. 
Practical algorithms for the Monge-Amp\`ere, Pucci, and $\sigma_2$ constraints can be found in \cite{dimitrios,pucci,glowinski} respectively. These algorithms are usually very inexpensive, as they are low-dimensional, converge in only a few iterations and can be parallelized.
We refer to \cite{caboussat,peruso2026,perusoestimates} for a detailed analysis of the least-squares formulation and the splitting algorithm.

Finally, following \cite{caboussat,pucci,lakkis}, we initialize $u^0$ as the solution to the Poisson equation: 
\begin{equation}\label{eq:initialization}
\Delta u^0 = 2\sqrt{f}\quad\text{in }\Omega,\quad u^0=\phi\quad \text{on }\partial\Omega.
\end{equation}
This provides a good initial guess when the eigenvalues of $D^2u$ are close to each other. Indeed, let $\lambda_1$ and $\lambda_2$ be the eigenvalues of $D^2u$. If $\lambda_1\approx\lambda_2$, then 
$$\Delta u=\text{\rm tr}(D^2u)\approx 2\lambda_1 \approx 2\sqrt{\det D^2u}=2\sqrt{f}.$$


\section{Extension to generalized Monge-Amp\`ere equations}\label{sec:ext}

In addition to the prototypical Dirichlet Monge-Amp\`ere problem and the Pucci's equation described earlier, we also consider two examples of generalized Monge-Amp\`ere equations: the prescribed Gauss curvature problem (also known as the Minkowski problem) and the optimal transport problem. 
In both examples the operator $F$ depends not only on $D^2u$ but also on $\nabla u$, and the optimal transport problem furthermore involves a different type of boundary conditions.

First, the prescribed Gauss curvature problem for the graph of a function $u:\Omega\subset\mathbb{R}^d\to\mathbb{R}$ consists of finding $u$ satisfying:
\begin{equation}\label{eq:prescribed-Gauss}
\begin{cases}
\det(D^2 u) = K(1 + |\nabla u|^2)^{(d+2)/2} \quad &\text{in } \Omega, \\
u = \phi \quad &\text{on } \partial \Omega.
\end{cases}
\end{equation}
where $K:\Omega\to(0,\infty)$ is a given, positive, function. 
Ellipticity of problem \eqref{eq:prescribed-Gauss} requires $u$ to be convex (so that $D^2u$ is positive semidefinite and the determinant is well-defined in the usual sense). 
To accommodate the dependence of $F$ on $\nabla u$ in our splitting algorithm, we modify the first substep \eqref{eq:firstmin} of our algorithm, and replace it by an explicit constraint on the current gradient as follows: 
 
\begin{equation}
\label{eq:nonlinearbis}
\mathbf{P}^n = \argmin_{\mathbf{Q} \in  L^2(\Omega; \mathbb{R}^{d\times d}_{\text{sym}}) }\left\{J(u^n,\mathbf{Q}), \quad\text{s.t. } F(\mathbf{Q},\nabla u^n) = 0\right\}. 
\end{equation}
The biharmonic step \eqref{eq:biharmonic} remains unchanged. 
An alternative to \eqref{eq:nonlinearbis} is to introduce an auxiliary variable $b$ approximating $\nabla u$ and solve an additional subproblem for $b$, see, e.g., \cite{minkowski}. 
In the numerical experiments presented in \Cref{sec:numres}, we have used the explicit-gradient formulation \eqref{eq:nonlinearbis} and numerical convergence was systematically observed. 

\begin{figure}[ht!]
\centering
\includegraphics[width=0.6\textwidth]{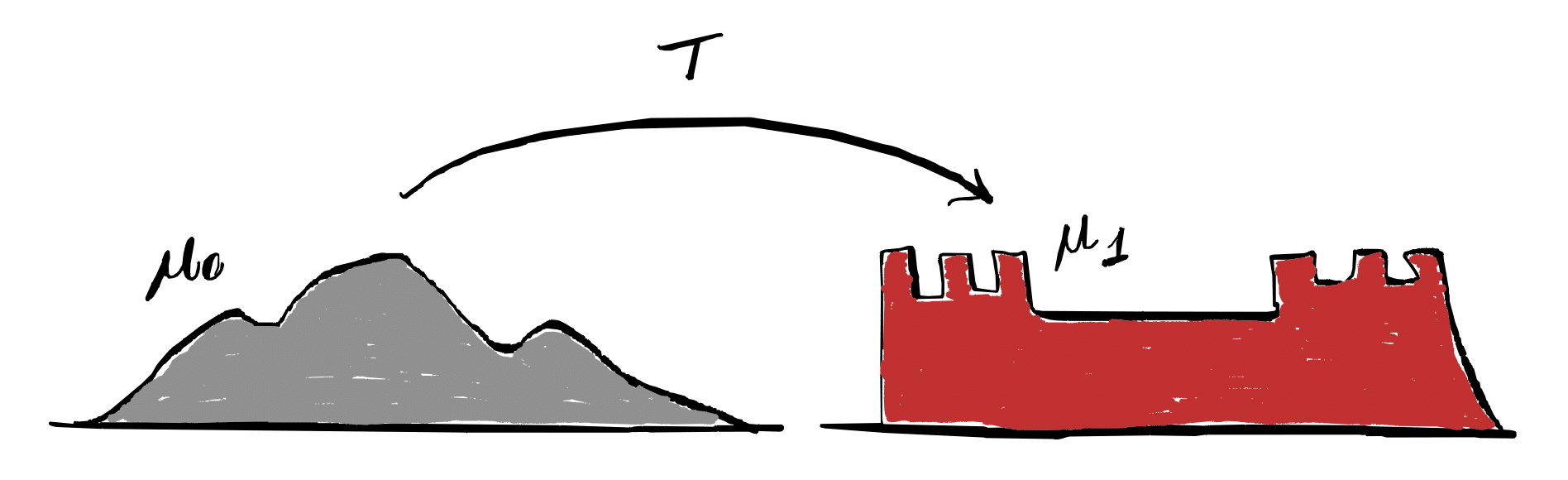}
\caption{Illustration of the Monge optimal transport problem (in 1D), which consists in finding the optimal transport map $T$ that transports the mass $\mu_0$ onto $\mu_1$.}	
\label{fig:OT}
\end{figure}

Second, the Monge-Amp\`ere optimal transport problem will be addressed. Originating in the work of Monge \cite{Monge1781}, the problem seeks an optimal way to transport mass between two given probability measures. 
More specifically, let $\mathcal{X}, \mathcal{Y}$ be two open bounded domains of $\mathbb{R}^d$. 
Given probability measures $\mu_0$ on $\mathcal{X}$ and $\mu_1$ on $\mathcal{Y}$, the objective is to find a measurable map $T:\mathcal{X}\to\mathcal{Y}$ whose pushforward 
matches $\mu_1$, i.e., 
\begin{equation}
\mu_1(A) = \mu_0 \left( T^{-1}(A) \right) \quad \forall \, A \subset \mathcal{Y} \text{ measurable},
\label{eq:pushforward}
\end{equation}
and that minimizes the transportation cost
\begin{equation}
\int_{\mathcal{X}} c(x, T(x)) \, d\mu_0(x) = \min_{S_{\#} \mu_0 = \mu_1} \int_{\mathcal{X}} c(x, S(x)) \, d\mu_0(x),
\label{eq:cost_minimization}
\end{equation}
where $c: \mathcal{X} \times \mathcal{Y} \rightarrow \mathbb{R}$ is the quadratic cost $c(x,y) = \frac{1}{2}|x-y|^2$. 
When \eqref{eq:pushforward} is satisfied, the map $T$ is called a transport map, and we use the notation $T_{\#}\mu_0=\mu_1$; if it also minimizes \eqref{eq:cost_minimization}, it is called an optimal transport map. 
\Cref{fig:OT} visualizes a transport map for $d=1$. 

When the cost is the quadratic cost, Brenier's theorem, see, {e.g.}, \cite{brenier}, gives a convenient PDE formulation of the optimal transport problem under standard regularity and support assumptions. If $\mu_0$ and $\mu_1$ are absolutely continuous probability measures on $\mathbb{R}^d$ with compact supports $\mathcal{X},\mathcal{Y}\subset\mathbb{R}^d$ and if the densities $\mu_0,\mu_1$ are bounded above and below on their supports, then:

\begin{enumerate}[i)]
  \item there exists a unique optimal transport map $T$ solving \eqref{eq:cost_minimization};
  \item there exists a convex potential $u:\mathcal{X}\to\mathbb{R}$ such that $T=\nabla u$ $\mu_0$-almost everywhere and $u$ satisfies the (generalized) Monge-Amp\`ere equation: 
  \begin{equation}
  \label{eq:Jacobian}
  \det(D^2 u)=\frac{\mu_0}{\mu_1(\nabla u)}\quad\text{for }\mu_0\text{-a.e. }x\in\mathcal{X},
  \end{equation}
  together with the {\em transport boundary condition} 
  \begin{equation}
  \label{eq:transportBC}
 \nabla u (\partial \mathcal{X}) = \partial \mathcal{Y}.
   \end{equation}
\end{enumerate}
As for the prescribed Gauss curvature problem, the operator $F$ appearing in the Monge-Amp\`ere equation depends on $\nabla u$, and the first substep \eqref{eq:firstmin} is replaced by the explicit-gradient strategy described in \eqref{eq:nonlinearbis}. 
Concerning boundary conditions \eqref{eq:transportBC}, numerical strategies have been proposed: for instance, one can project $\nabla u$ onto the set of feasible maps as an extra minimization subproblem \cite{prins}. In this work we follow the approach detailed in \Cref{sub:ot}.


\section{The Deep Ritz method}\label{sec:deepritz}
In this section, we detail how to solve \eqref{eq:biharmonic} by using a neural network approximation. 
Owing to the variational structure of the problem, we adopt the Deep Ritz method \cite{deepritz}. 
We first briefly introduce the basic framework of neural networks and then present the Deep Ritz method together with a novel adaptive sampling technique in detail. 
Finally, we explain how to extend the method from Dirichlet to transport boundary conditions.


\subsection{Neural Networks and Input Convex Neural Networks}\label{sub:nn}

We consider a feedforward neural network (NN) with one input layer, $L \geq 1$ hidden layers, and one output layer. 
The width of the $j$th layer is denoted by $N_j$ for $j = 0,\dots,L+1$, where $N_0$ and $N_{L+1}$ specify the input and output dimensions. 
Each layer is parametrized by a weight matrix $W^{(l)} \in \mathbb{R}^{N_{l+1} \times N_l}$ and a bias vector $b^{(l)} \in \mathbb{R}^{N_{l+1}}$ for $l = 0,\dots,L$.  
Given an input $x$, the forward propagation is defined by
\begin{align*}
    x^0 &= x,\\
    x^l &= \sigma(W^{(l-1)} x^{l-1} + b^{(l-1)}), \quad 1 \leq l \leq L,\\
    y &= W^{(L)} x^L + b^{(L)},
\end{align*}  
where $\sigma$ is a nonlinear activation function applied componentwise. 
Collecting all parameters into $\theta := \{W^{(l)}, b^{(l)}\}_{l=0}^L$, the network defines a approximation $u_{\text{NN}}(x;\theta)$ of the solution $u$.   

We denote by $\mathcal{V}^\sigma_{\{N_j\}}$ the class of neural networks with activation function $\sigma$ and architecture specified by the widths $\{N_j\}_{j=0}^{L+1}$. 

\medskip

Solutions of some fully nonlinear PDEs must satisfy structural constraints; for example, we enforce the solutions of the Monge-Amp\`ere equation to be convex. 
To enforce such convexity in a neural network representation, one may employ Input Convex Neural Networks (ICNNs) \cite{icnn}, which achieve convexity by incorporating specific architectural restrictions and parameter constraints. %
An ICNN with input $x$ and $L$ hidden layers is defined by
\begin{align*}
x^0 &= x, \\
x^1 &= \sigma(L^{(0)}x^0 + b^{(0)}), \\
x^l &= \sigma(W^{(l-1)}x^{l-1} + L^{(l-1)}x^0 + b^{(l-1)}), \quad 2 \leq l \leq L, \\
y &= W^{(L)}x^{L} + L^{(L)}x^0 + b^{(L)},
\end{align*}
where $\sigma$ is an increasing convex activation function. The parameters are constrained as follows:  
$W^{(l)} \in \mathbb{R}^{N_{l+1}\times N_l}$ have nonnegative entries, ensuring that each hidden layer produces a convex combination of its inputs; 
$L^{(l)} \in \mathbb{R}^{N_{l+1}\times N_0}$ are unconstrained skip-connection matrices from the input;
and $b^{(l)} \in \mathbb{R}^{N_{l+1}}$ are bias terms.  
\Cref{fig:icnn} visualizes such an ICNN.

\begin{figure}[ht!]
\centering
\resizebox{0.9\linewidth}{!}{%
\begin{tikzpicture}[>=Stealth, node distance=8mm, every node/.style={scale=0.85}]
  \node[draw=black, fill = blue!20, rectangle, rounded corners, align=center] (x) {Input \\ $x$};

  \node[draw=black, fill = pink!20, rectangle, rounded corners, right=15mm of x, align=center] (h1)
        {$\sigma(L^{(0)}x + b^{(0)})$};

  \node[draw=black, fill = red!20, rectangle, rounded corners, right=18mm of h1, align=center] (h2)
        {$\sigma(W^{(1)}x^{(1)} + L^{(1)}x + b^{(1)})$};

  \node[draw=black, fill = red!20, rectangle, rounded corners, right=18mm of h2, align=center] (hL)
        {$\cdots$};

  \node[draw=black, fill = green!20, rectangle, rounded corners, right=15mm of hL, align=center] (y) {Output \\ $y$};

  \draw[->, dashed] (x.east) to[bend left=10] node[above, midway]{$L^{(0)}$} (h1.west);
  \draw[->, thick] (h1.east) to[bend left=10] node[above, midway]{$W^{(1)}\!\ge 0$} (h2.west);
  \draw[->, thick] (h2.east) to[bend left=10] node[above, midway]{$W^{(2)}\!\ge 0$} (hL.west);
  \draw[->, thick] (hL.east) to[bend left=10] node[above, midway]{$W^{(L)}\!\ge 0$} (y.west);

  \draw[->, dashed] (x.east) to[bend right=20] node[below, midway]{$L^{(1)}$} (h2.west);
  \draw[->, dashed] (x.east) to[bend right=26] node[below, midway]{$L^{(2)}$} (hL.west);
  \draw[->, dashed] (x.east) to[bend right=32] node[below, midway]{$L^{(L)}$} (y.west);
\end{tikzpicture}%
}
\caption{Schematic of an Input Convex Neural Network (ICNN) and notation.}
\label{fig:icnn}
\end{figure}
The convexity of the solution is then guaranteed by two elementary facts:  
\begin{enumerate}[i)]
  \item nonnegative linear combinations of convex functions are convex;  
  \item the composition of a convex function with an increasing convex function yields a convex function.  
\end{enumerate}
A specific aspect of ICNNs is the use of \textit{passthrough} layers $\{L^{(l)}\}_{l}$, which provide direct connections from the input $x$ to hidden units in deeper layers. 
In standard feedforward networks, such connections are unnecessary, since hidden units can pass information forward without restriction. 
In contrast, ICNNs impose a non-negativity constraint on the weights $W^{(l)}$, thereby limiting the ability of hidden units to reproduce identity mappings. 
To address this limitation, passthrough layers are introduced explicitly, and they play a key role in ensuring the network’s representational power \cite{Chen2018}.


\subsection{The Deep Ritz method for the variational problem \eqref{eq:biharmonic}}\label{sec:deepritz_biharm}

The variational subproblem \eqref{eq:biharmonic} is linear and reads:
$$u^{n+1} = \argmin_{\substack{v\in H^2(\Omega)\\ v|_{\partial\Omega}=\phi}}\int_{\Omega}|D^2v - \mathbf{P}^n|^2.$$
The problem is well-posed since $\int_\Omega |D^2(\cdot)|^2$ defines a norm in $H^2(\Omega)\cap H^1_0(\Omega)$.

The Deep Ritz method \cite{deepritz} employs neural networks to approximate the solutions of PDEs via their variational formulation, thereby avoiding traditional mesh-based discretizations. 
 The preference for the Deep Ritz method over other physics-guided approaches, such as PINNs \cite{karniadakis2021}, is motivated by the natural variational structure of the problem, and the facts that no higher-order derivatives need to be computed, and no additional regularity assumptions are imposed. 
The solution to \eqref{eq:biharmonic} is thus approximated by a deep neural network solution $u^{n+1}_{NN}$, which satisfies
\begin{equation}  
\label{eq:losscont}
u^{n+1}_{NN} = \argmin_{v\in  \mathcal{V}^\sigma_{\{N_j\}}} \left[ \int_\Omega |D^2v - \mathbf{P}^n|^2 +\lambda\int_{\partial\Omega}(v -\phi)^2 \right].  
\end{equation}  
Since, in general, $\mathcal{V}^\sigma_{\{N_j\}} \not\subset H^2(\Omega) \cap H^1_\phi(\Omega)$, the penalty term in \eqref{eq:losscont} allows to enforce weakly the Dirichlet boundary condition (see \cite{essential_DR} for an alternative treatment of essential boundary conditions).
However, for sufficiently smooth activation functions $\sigma$, we have $\mathcal{V}^\sigma_{\{N_j\}} \subset H^2(\Omega)$, which provides the network with some degree of regularity. 

In practice, the {\em loss function} is obtained by Monte Carlo discretization of \eqref{eq:losscont}. 
Given a set of $n_c$ uniformly sampled interior collocation points $\{x_{c,i}\}_{i=1}^{n_c} \subset \Omega$ and $n_b$ uniformly sampled boundary points $\{x_{b,i}\}_{i=1}^{n_b} \subset \partial\Omega$, the continuous loss function is approximated by
$$
\mathcal{L}(v) 
= \underbrace{\frac{1}{n_c} \sum_{i=1}^{n_c}  \left| D^2 v(x_{c,i}) - \mathbf{P}^n(x_{c,i}) \right|^2}_{:=\mathcal{L}_{PDE}(v)}
+ \lambda \underbrace{\frac{1}{n_b} \sum_{i=1}^{n_b} \left( v(x_{b,i}) - \phi(x_{b,i}) \right)^2}_{:=\mathcal{L}_{BC}(v)}.
$$
The optimization problem then becomes: 
$$
u^{n+1}_{\mathrm{NN}} \approx \argmin_{v \in \mathcal{V}^\sigma_{\{N_j\}}} \mathcal{L}(v).
$$
A possible drawback of the proposed method is that the least-squares Deep Ritz formulation requires solving many separate and sequential optimization problems (one per splitting iteration), whereas a PINN approach approximates the full Monge-Amp\`ere equation directly and might appear computationally cheaper and simpler to implement \cite{pinnsconvex}. 
In practice, however, we expect that the Deep Ritz training requires only a small number of epochs after the first few splitting iterations, due to the continuation effect provided as the solution at the previous splitting iteration should give a good initialization for the next Deep Ritz iteration.

Consequently, the number of epochs used for each Deep Ritz solve is reduced as the iterations progress (see \Cref{sec:impl} for implementation details). 

\subsection{Adaptive sampling}\label{sub:adaptivesampling}
The accuracy of a Deep Ritz approximation depends not only on the approximation capability of the neural network and on the optimization error, but also on the accuracy with which the variational functional is approximated by Monte Carlo quadrature. Uniform sampling can be inefficient when the loss density is strongly localized. Adaptive sampling strategies have therefore been investigated extensively for physics-informed neural networks and related neural PDE solvers; see, for example, \cite{adaptive-comp} and the references therein. In the context of the Deep Ritz method, \cite{adaptive_deep_ritz} approximate the energy integrand by a generative density model and use the resulting density for importance sampling, while \cite{adaptive_deep_ritz_elasticity} employ a residual-based adaptive quadrature method for linear elasticity.

The strategy considered here follows the same general importance-sampling principle as \cite{adaptive_deep_ritz}. Related to methods based on an auxiliary generative model or on an adaptively refined integration mesh, our implementation approximates the sampling density by a piecewise-constant function obtained from a small number of seed evaluations, following \cite{efficient_sampling}. This produces a simple mesh-free procedure with a relatively small density-update cost. Its limitation is that the resulting approximation of the ideal importance density may be coarse; the dependence of the method on the number of seeds is investigated in \Cref{sec:numres}.

Given any strictly positive probability density \(q\) on \(\Omega\), the Lebesgue integral defining the PDE loss admits the exact importance-sampling representation
\[
\int_\Omega |D^2v(x)-\mathbf P^n(x)|^2\,dx
=
\mathbb E_{X\sim q}
\left[
\frac{|D^2v(X)-\mathbf P^n(X)|^2}{q(X)}
\right].
\]

Consequently, for independent samples \(X_i\sim q\), the Monte Carlo estimator
\begin{equation}
\widehat{\mathcal L}_{\mathrm{PDE}}(v)
=
\frac{1}{n_c}
\sum_{i=1}^{n_c}
\frac{
|D^2v(X_i)-\mathbf P^n(X_i)|^2
}{
q(X_i)
},
\label{eq:IS-loss}
\end{equation}
is unbiased. The density minimizing the variance of this estimator, in the formal sense, is given by
\[
q^\star(x)
=
\frac{d(x)}
{\displaystyle\int_\Omega d(y)\,dy},
\qquad
d(x):=|D^2v(x)-\mathbf P^n(x)|^2.
\]
Here, \(v\) denotes the current network approximation at the time the sampling distribution is updated.

Directly evaluating \(d(x)\) over all collocation points at every iteration, however, can be computationally prohibitive. Motivated by the optimal importance density $q^\star$, we construct an adaptive stratification of the domain. The purpose is to allocate more collocation points to regions where the current loss density is large, while a quadrature correction ensures that the original energy loss is correctly approximated. Let \(\{X_s\}_{s=1}^S\) be a small set of seed points, with \(S\ll n_c\). The seeds induce a Voronoi tessellation of \(\Omega\) through a nearest-neighbor assignment. Denoting the corresponding Voronoi cells by \(C_s\), we approximate the loss density within each cell by its value at the associated seed, \(d(X_s)\). In particular, we define the number of points sampled on each cell by 
$$n_s = n_c \frac{d(X_s)}{\sum_{j=1}^S d(X_j)}$$
and let $|C_s|$ be the volume of the cell, then we approximate the loss \eqref{eq:IS-loss} by 
$$
\widetilde{\mathcal L}_{\mathrm{PDE}}(v)
=\sum_{s = 1}^S \frac{|C_s|}{n_s}\sum_{i=1}^{n_s}
|D^2v(X_{s,i})-\mathbf P^n(X_{s,i})|^2, \qquad X_{s,i}\sim \mathcal{U}(C_s).
$$
The cell volumes \(|C_s|\) are estimated using Monte Carlo samples drawn from an initial uniform distribution over the collocation domain. The approach is visualized in \Cref{sec:numres}, for example in \Cref{fig:exp_resampled}.

Compared with evaluating the indicator at all \(n_c \) points, this approach requires only \(S\) Hessian evaluations for each density update, in addition to the nearest-neighbor assignment. Its effectiveness depends on how accurately the seed-induced tessellation resolves the spatial variation of the loss density; this dependence is examined numerically in \Cref{sec:numres}.

\subsection{Extension to optimal transport problems}\label{sub:ot}
Within the Deep Ritz framework, a modification of the loss function $\mathcal{L}(v)$ allows to weakly handle alternative boundary conditions, such as, e.g., those arising in the optimal transport problem introduced in \Cref{sec:ext}, following the method proposed in \cite{pinnsconvex}. 
Instead of enforcing $u|_{\partial\Omega}=\phi$ via a penalty term, we must now impose the transport boundary condition \eqref{eq:transportBC}.
Note that
$$
\nabla u(\partial \mathcal{X}) = \partial \mathcal{Y} \quad \Longleftrightarrow \quad d_H(\nabla u(\partial \mathcal{X}), \partial \mathcal{Y}) = 0,
$$
where $d_H$ denotes the Hausdorff distance (with respect to the Euclidean metric) defined by 
$$
d_H(A, B) := \max \Big\{ \sup_{a \in A} \mathrm{dist}(a, B), \;\; \sup_{b \in B} \mathrm{dist}(b, A) \Big\}, \qquad A,B \subset \mathbb{R}^n,
$$
and
$$
\mathrm{dist}(a, B) := \min_{b \in B} \|a - b\|_2.
$$
The two components of the Hausdorff distance enforce the two set inclusions
$$\nabla u(\partial\mathcal X)\subseteq\partial\mathcal Y\qquad\text{and}\qquad\partial\mathcal Y\subseteq\nabla u(\partial\mathcal X),$$ 
respectively. Hence, their simultaneous vanishing enforces $\nabla u(\partial\mathcal X)=\partial\mathcal Y$.
Thus the continuous constrained problem can be written as
$$
\min_{v \in \mathcal{V}^\sigma_{\{N_j\}}} \left\{ \int_\Omega |D^2 v - \mathbf{P}^n|^2 \;\; \text{s.t. } d_H(\nabla v(\partial \mathcal{X}), \partial \mathcal{Y}) = 0 \right\}.
$$
We approximate the Hausdorff distance in a discrete setting as follows. 
Let $\{x_{\mathcal{X},i}\}_{i=1}^{n_{b\mathcal{X}}}\subset\partial\mathcal{X}$ and $\{y_{\mathcal{Y},i}\}_{i=1}^{n_{b\mathcal{Y}}}\subset\partial\mathcal{Y}$ be finite sets of boundary points of cardinalities $n_{b\mathcal{X}}$ and $n_{b\mathcal{Y}}$, respectively. We then define the discrete loss function by 
$$
\mathcal{L}(v) = \mathcal{L}_{{PDE}}(v) + \lambda \mathcal{L}_{{OT}}(v),
$$
where the PDE term $\mathcal{L}_{{PDE}}(v)$ remains unchanged and the transport-boundary term is now defined as: 
\begin{equation}
\mathcal{L}_{{OT}}(v) =
\frac{1}{n_{b\mathcal{X}}} \sum_{i=1}^{n_{b\mathcal{X}}} 
\mathrm{dist}\Big(\nabla v(x_{\mathcal{X},i}), \{y_j\}_{j=1}^{n_{b\mathcal{Y}}}\Big)^2
+ \frac{1}{n_{b\mathcal{Y}}} \sum_{i=1}^{n_{b\mathcal{Y}}} 
\mathrm{dist}\Big(\nabla v(\{x_j\}_{j=1}^{n_{b\mathcal{X}}}), y_{\mathcal{Y},i}\Big)^2,
\label{eq:tildeEb}
\end{equation}
The first term in \eqref{eq:tildeEb} penalizes deviations from injectivity of the map $\nabla v$ on the source boundary, while the second term enforces surjectivity onto the target boundary.

The Hausdorff distance formulation is general and can be applied in a variety of settings, although it can be computationally demanding. 
In some special cases the transport boundary condition admits a simpler caracterization; for instance, when the source and target are identical square domains, the optimal-transport condition is equivalent to a Neumann boundary condition. In such situations the boundary terms can be simplified and the solver can be optimized. 
We do not pursue this direction in the present work.


\section{Implementation details}\label{sec:impl}
The splitting algorithm and the Deep Ritz solver have been implemented in \texttt{PyTorch}. 
When considering the Monge-Amp\`ere equation, we approximate $u^0$ in \eqref{eq:initialization} by a neural network trained with the Deep Ritz method. The same initialization problem \eqref{eq:initialization} is also applied for the Pucci's case \eqref{intro_pucci} (which corresponds to the case $f=0$ and $\alpha=1$). For the Minkowski case \eqref{eq:prescribed-Gauss}, we set $\nabla u^0 = \textrm{Id}$ on the right-hand side before solving the same problem \eqref{eq:initialization}. For the OT problem \eqref{eq:Jacobian}-\eqref{eq:transportBC}, we initialize the problem such that $\nabla u^0 = \textrm{Id}$.

At each outer iteration, the minimization problem \eqref{eq:firstmin} is solved pointwise at the collocation points, with existing local nonlinear solvers described in \cite{dimitrios,pucci,glowinski} respectively. Since these problems are low-dimensional and completely decoupled across the collocation points, this step has a negligible computational cost compared with the neural network optimization. In our timing experiments, the pointwise minimization accounts for at most approximately $0.12\%$ of the CPU time required for one variational update. The biharmonic problem \eqref{eq:biharmonic} is approximated with the Deep Ritz method and fully connected feed-forward networks with four hidden layers and ten neurons per layer, unless otherwise stated. For problems requiring convex solutions, such as, e.g., the Monge-Amp\`ere equation, we employ ICNNs. Following \cite{NYSTROM,pinnsconvex}, the ICNNs use the softplus activation function $\sigma(x) = \log(1+e^{x})$, which is smooth, increasing, and strictly convex. 
Convexity is enforced by constraining the relevant weight matrices to be nonnegative. 
To guarantee the positivity of $W^{(l)} \in \mathbb{R}^{N_{l+1} \times N_{l}}$, the weights are initialized by squaring them element-wise \cite{pinnsconvex,NYSTROM}.

\medskip

To improve training stability and accuracy, we adopt the common two-stage optimizer schedule used in physics-informed and Deep Ritz methods \cite{pinnstraining}. 
Each solver is first pre-trained with the Adam optimizer to achieve rapid initial progress, and then refined with the L-BFGS optimizer for improved convergence. 
During the outer splitting iterations, the L-BFGS number of epochs is gradually reduced: the first iterations employ up to $70$ L-BFGS epochs, which are progressively decreased to as few as $4$ epochs in later iterations. 
\Cref{tab:training} summarizes the training configuration used in our experiments. Finally, in the following sections, error plots are averaged on $10$ independent simulations and the shadowed area represents the 5th-95th quantile.

\begin{table}[ht!]
  \centering
  \caption{Training configuration for Deep Ritz solvers.}
  \label{tab:training}
  \begin{tabular}{l|l}
    Main optimizer                  & L-BFGS \\
    Line search                     & strong Wolfe \\
    L-BFGS history size             & 25\\
    Max iterations / call          & 20 \\
    Gradient tolerance (L-BFGS)     & $10^{-7}$ \\
    Batch size                      & full-batch 
\end{tabular}
\end{table}


\section{Numerical Results for Dirichlet problems}\label{sec:numres}

\subsection{Monge-Amp\`ere equation in 2D : preliminary test}
\begin{figure}[ht!]
\centering
\begin{subfigure}{0.40\linewidth}
\centering
\includegraphics[width=0.8\linewidth]{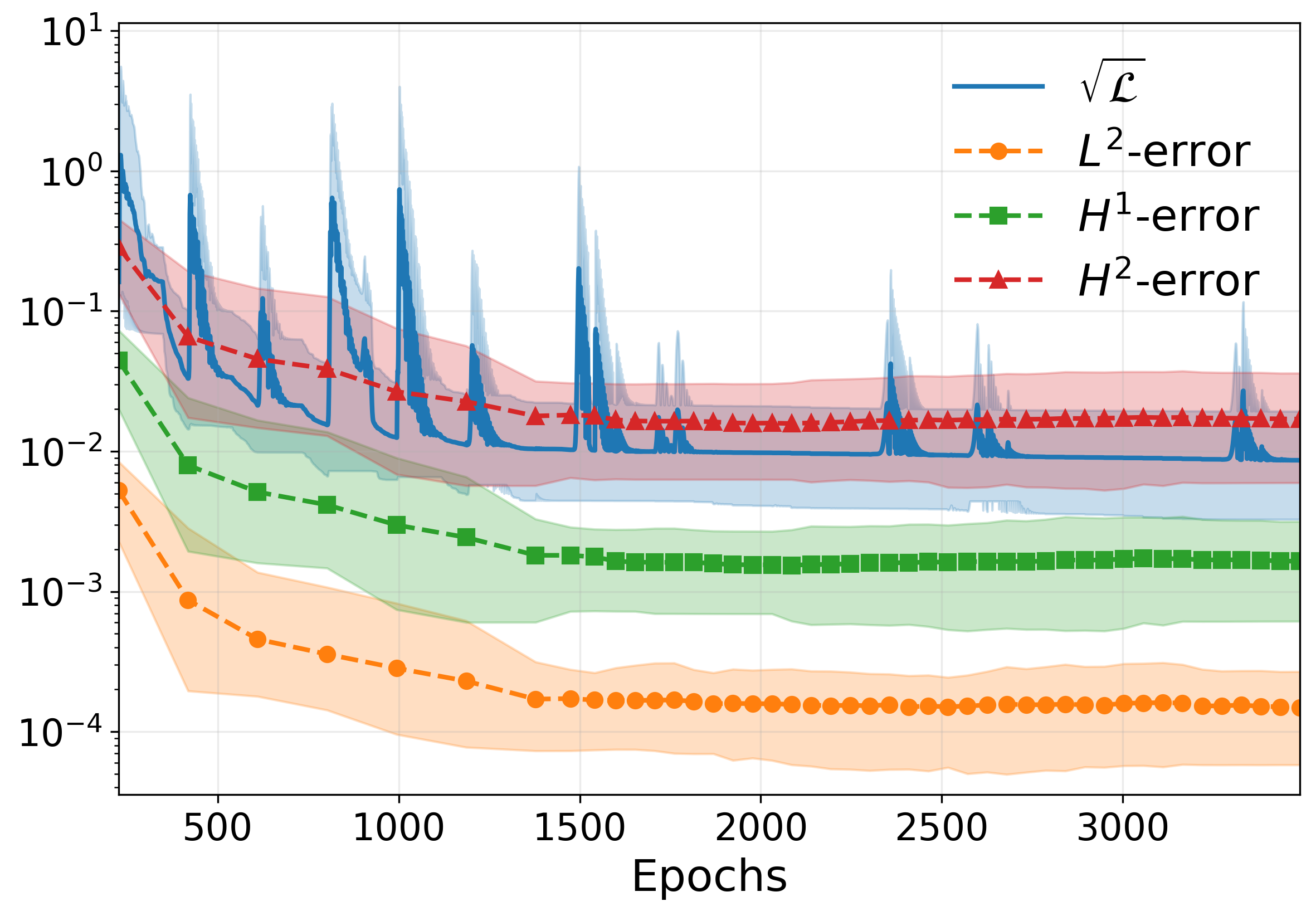}
\caption{PINNs.}
\end{subfigure}
\begin{subfigure}{0.40\linewidth}
\centering
\includegraphics[width=0.8\linewidth]{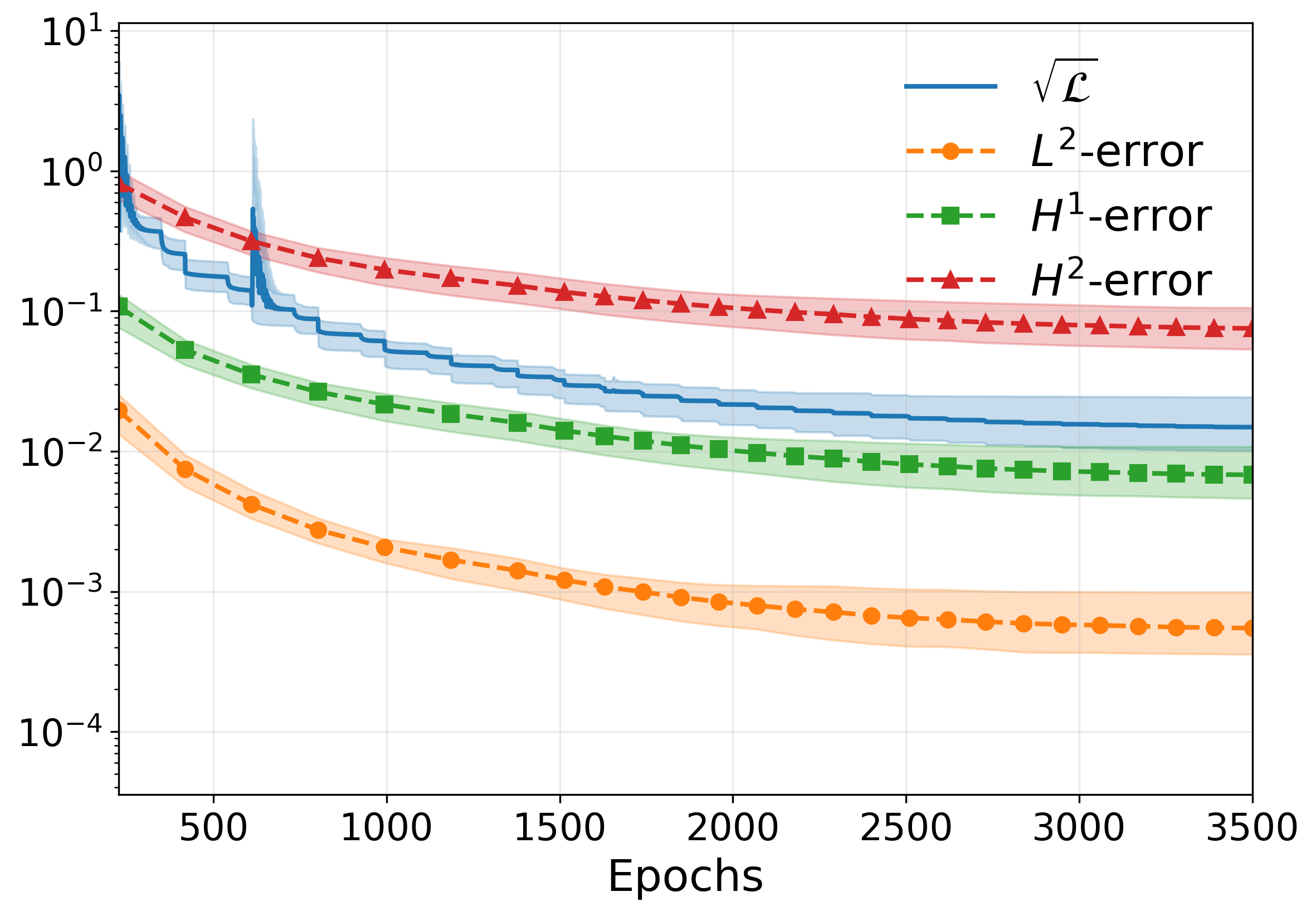}
\caption{Deep Ritz method.}
\end{subfigure}
\caption{First exponential test case, with $\alpha =1$ and 3000 collocation points. 
We illustrate the loss function for each epoch and the errors between $u$ and $u_{NN}$ computed only at the end of each splitting iteration. 
The same total number of epochs is used for PINNs.}
\label{fig:exp_1_comp}
\end{figure}

\begin{figure}[ht!]
\centering
\begin{subfigure}{0.32\linewidth}
\centering
\includegraphics[width=0.9\linewidth]{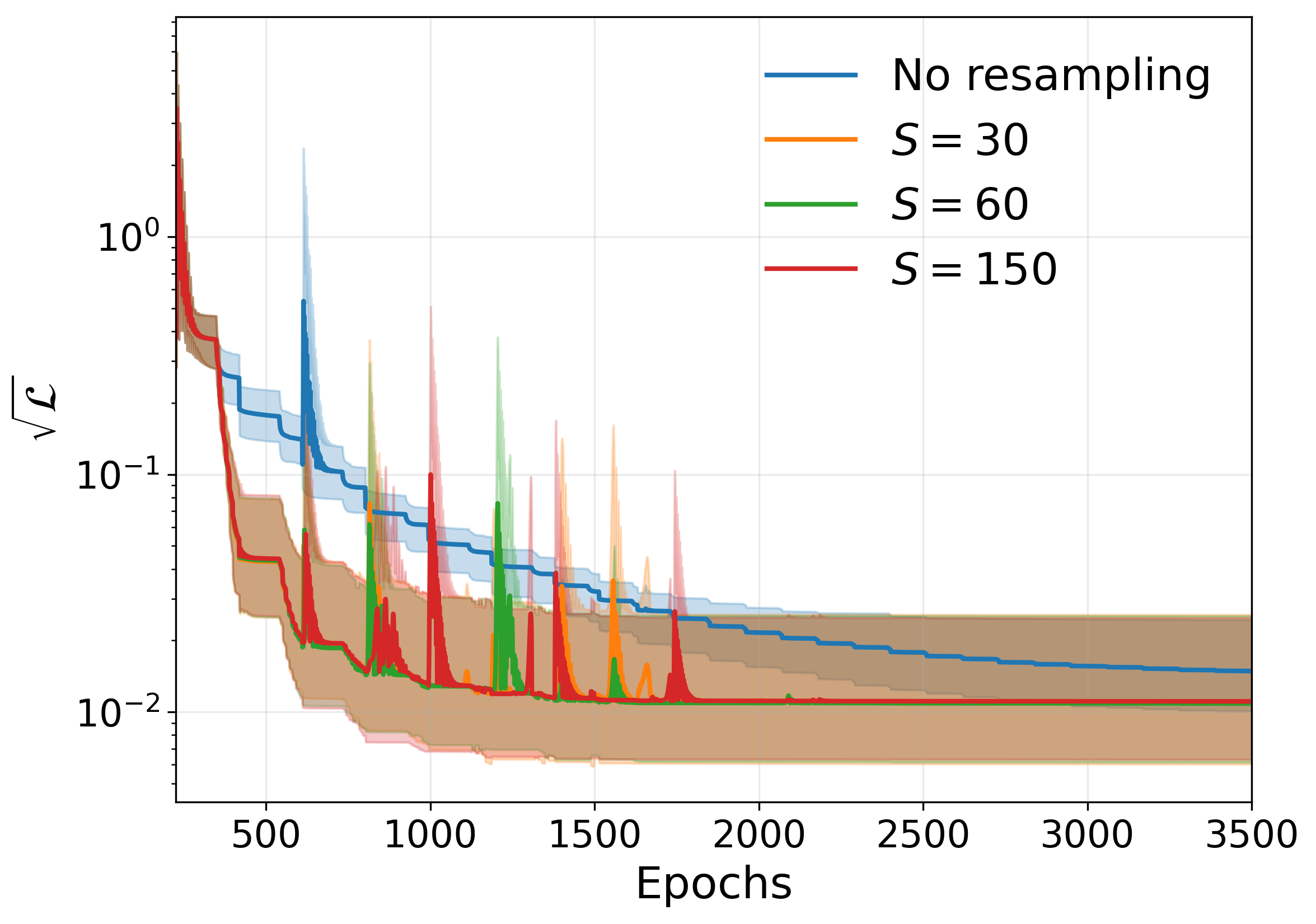}
\caption{Loss function.}
\end{subfigure}
\begin{subfigure}{0.32\linewidth}
\centering
\includegraphics[width=0.9\linewidth]{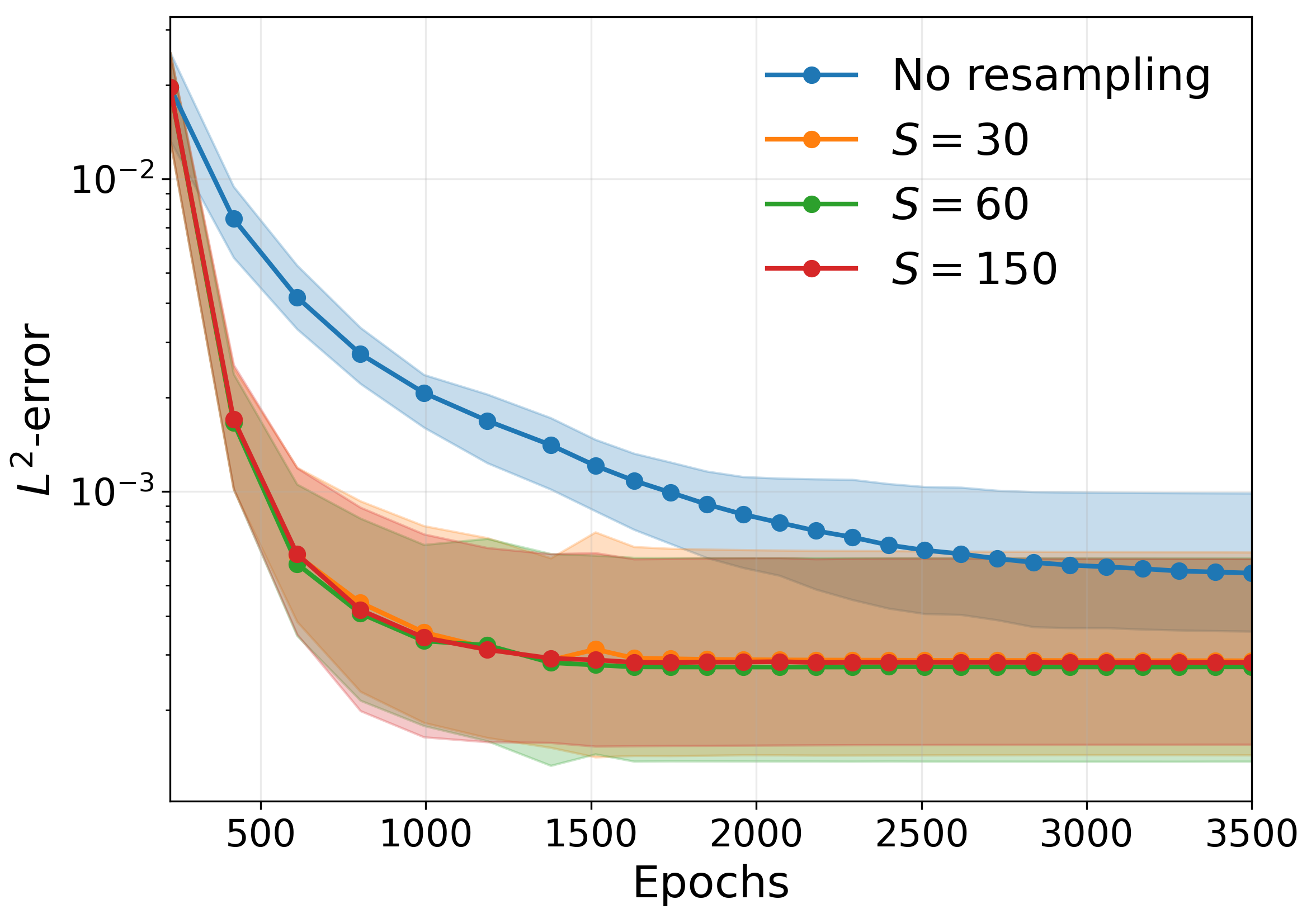}
\caption{$L^2$ error.}
\end{subfigure}
\begin{subfigure}{0.32\linewidth}
\centering
\includegraphics[width=0.9\linewidth]{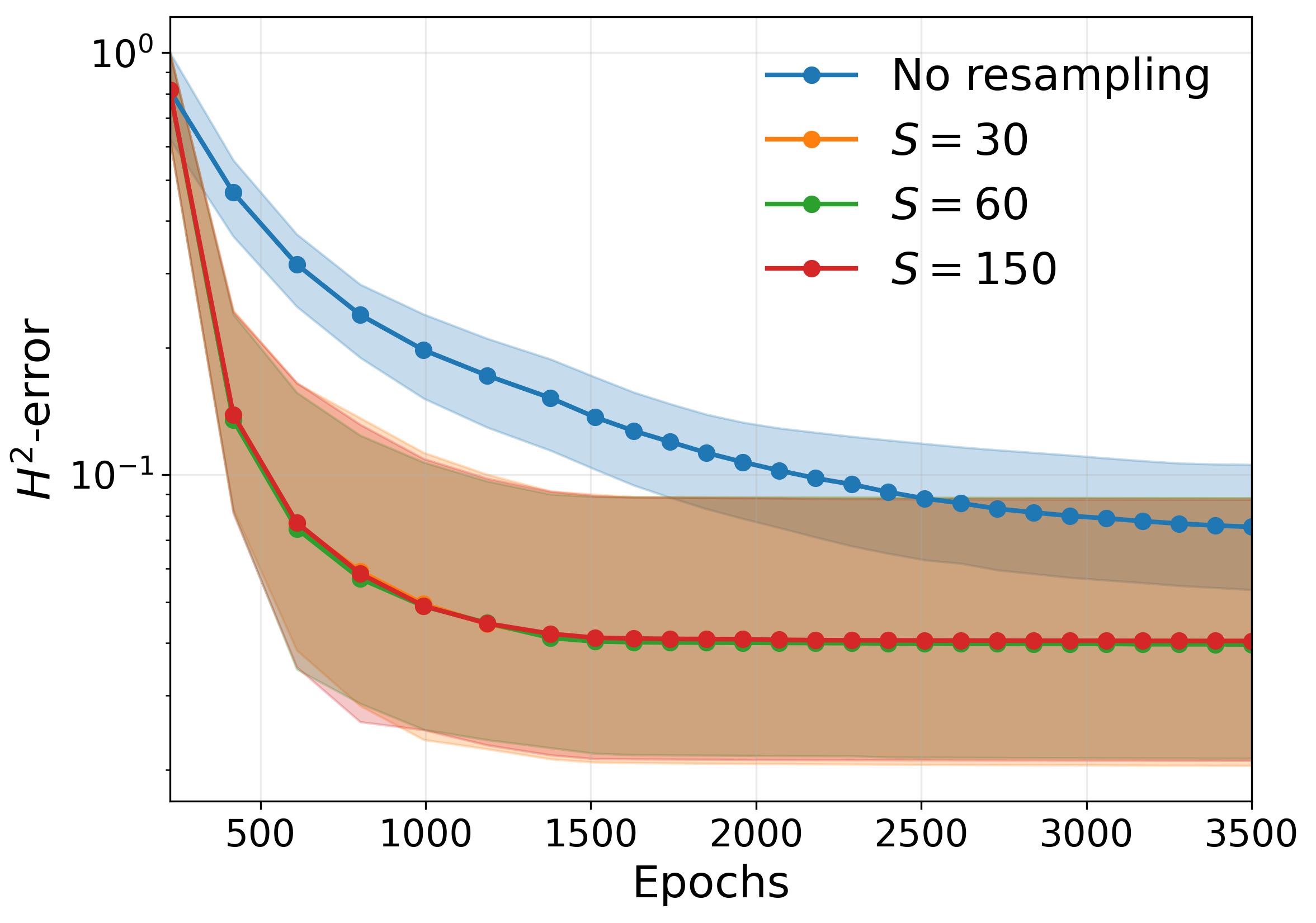}
\caption{$H^2$ error}
\end{subfigure}
\caption{First exponential test case, with $\alpha =1$. Deep Ritz method. Loss function (left), $L^2$ error (center) and $H^2$ error (right) for various values of the seed $S$. }
\label{fig:exp_adapt}
\end{figure}

\begin{figure}[ht!]
\centering
\begin{subfigure}{0.95\linewidth}
\centering
\includegraphics[width=0.95\linewidth]{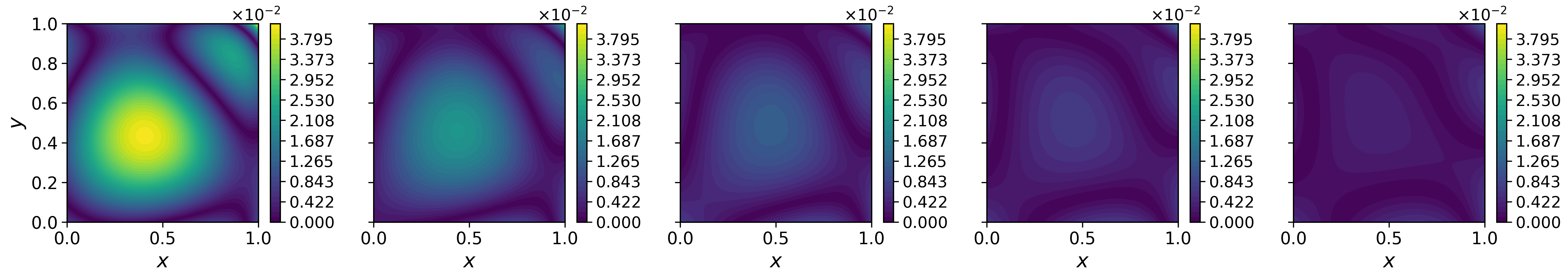}
\caption{No re-sampling.}
\end{subfigure}
\vspace{0.5cm}
\begin{subfigure}{0.95\linewidth}
\centering
\includegraphics[width=0.95\linewidth]{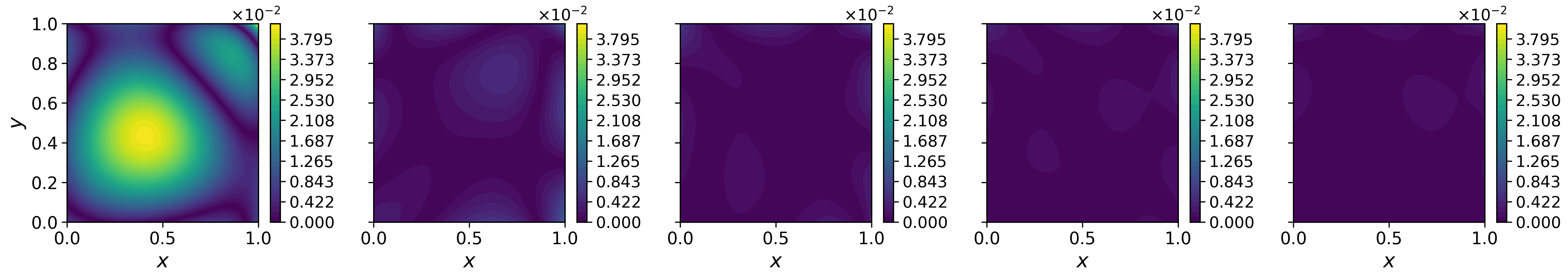}
\caption{Adaptive sampling.}
\end{subfigure}
\caption{First exponential test case, with $\alpha =1$, 3000 collocation points and 300 boundary points. 
Pointwise absolute error at the end of each splitting iteration.
Top row: no adaptive sampling; Bottom row: with adaptive sampling. }
\label{fig:exp_pointwise}
\end{figure}

\begin{figure}[ht!]
\centering
\includegraphics[width=0.9\linewidth]{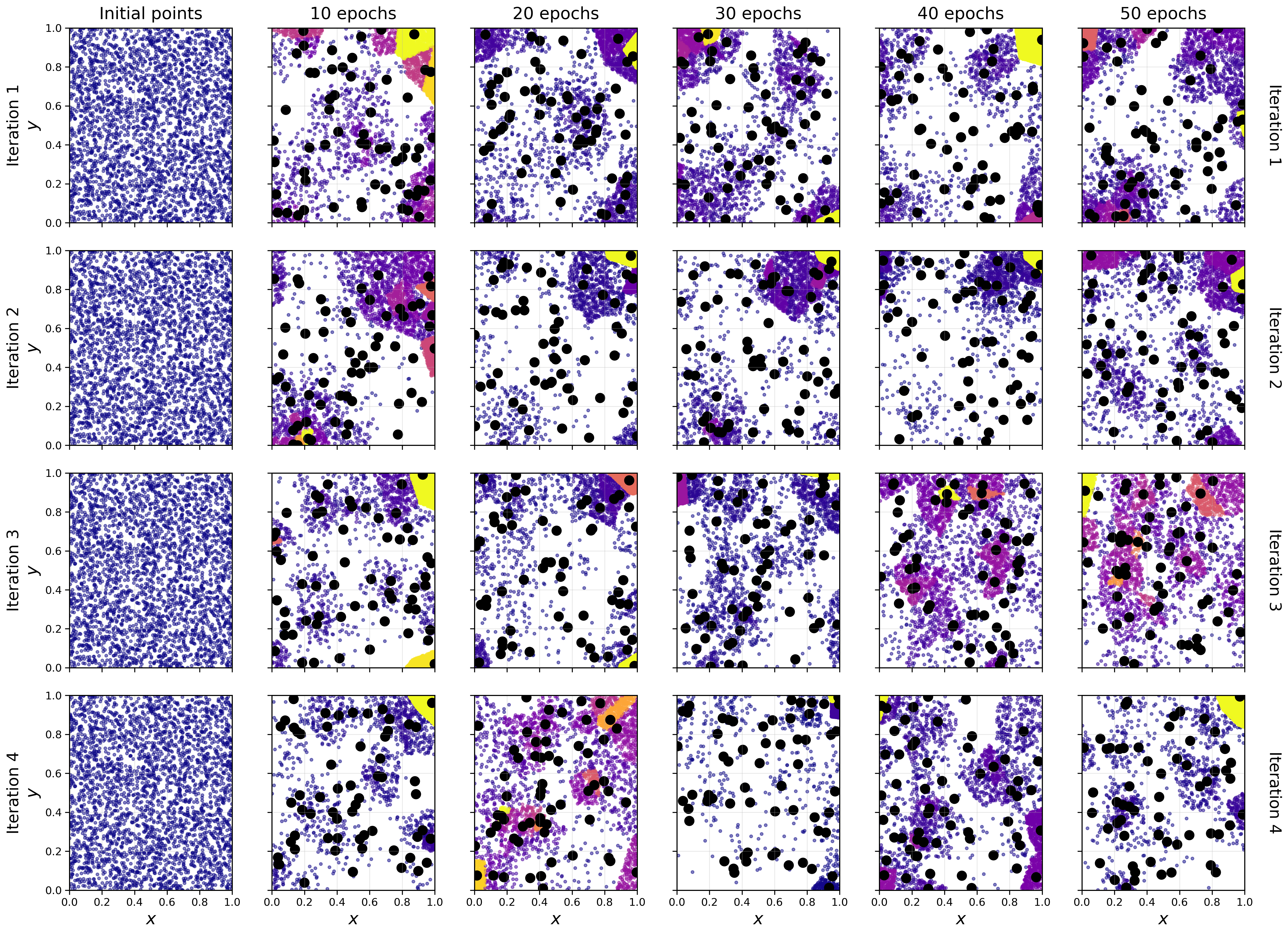}
\caption{First exponential test case, with $\alpha =1$, 3000 collocation points and 300 boundary points. 
Distribution of the resampled points (the seeds are illustrated in black).}
\label{fig:exp_resampled}
\end{figure}

Let $\Omega = [0,1]^2$. We consider the Monge-Amp\`ere problem $\det(D^2u)=f$ with Dirichlet boundary conditions $u = \phi$. 
The data $f$ and $\phi$ are such that the exact solution is given by
$$
u(x, y) = e^{(\alpha/2)(x^2 + y^2)},\quad\alpha>0.
$$
The case $\alpha=1$ is a standard benchmark. We compare the Deep Ritz approach with the PINN method of \cite{NYSTROM}, using $3000$ collocation points and $\lambda=100$. As shown in \Cref{fig:exp_1_comp}, PINNs converges faster and Deep Ritz shows smaller variability. Adaptive sampling further improves the accuracy and stability of the Deep Ritz approximation; increasing the number of seeds beyond moderate values provides little additional benefit (\Cref{fig:exp_adapt}). The corresponding pointwise errors and the evolution of the collocation points are shown in \Cref{fig:exp_pointwise,fig:exp_resampled}.

\begin{figure}[ht!]
\centering
\begin{subfigure}{0.40\linewidth}
\centering
\includegraphics[width=0.8\linewidth]{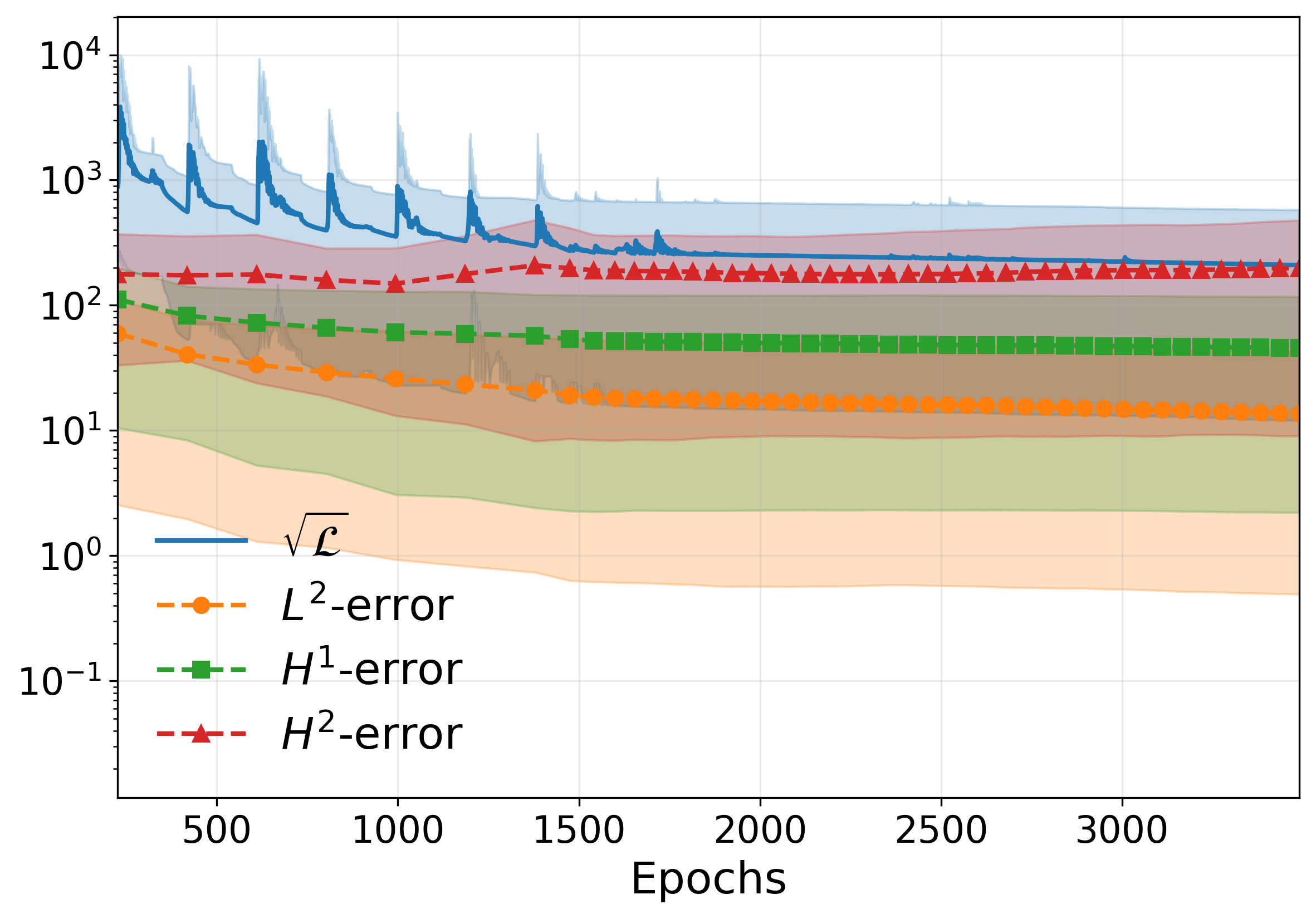}
\caption{PINNs.}
\end{subfigure}
\begin{subfigure}{0.40\linewidth}
\centering
\includegraphics[width=0.8\linewidth]{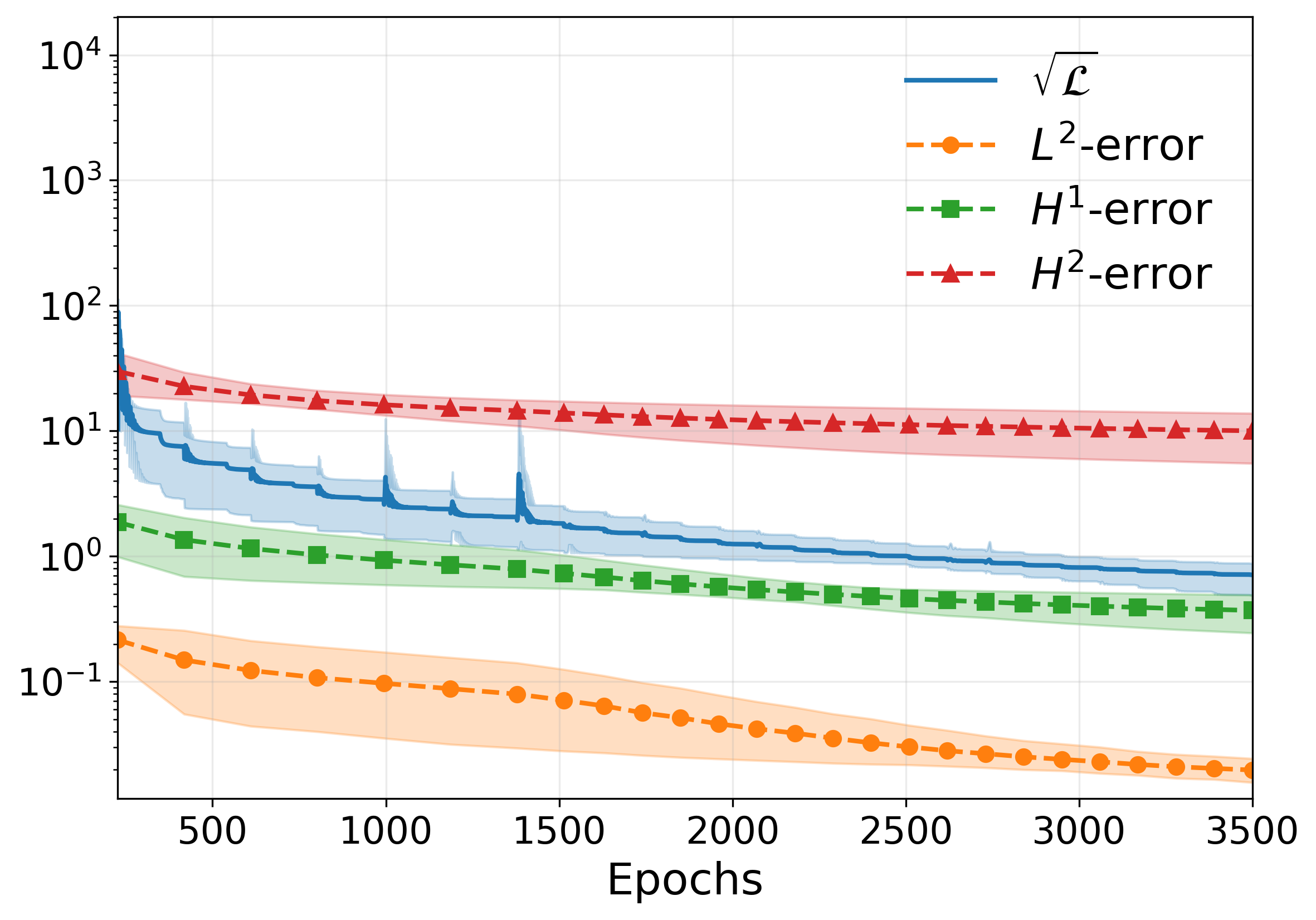}
\caption{Deep Ritz method.}
\end{subfigure}
\caption{First exponential test case, with $\alpha =4$ and 3000 collocation points. 
We illustrate the loss function for each epoch and the errors between $u$ and $u_{NN}$ computed only at the end of each splitting iteration. 
The same total number of epochs is used for PINNs. 
}
\label{fig:exp_4_comp}
\end{figure}

\begin{figure}[ht!]
\centering
\begin{subfigure}{0.32\linewidth}
\centering
\includegraphics[width=0.9\linewidth]{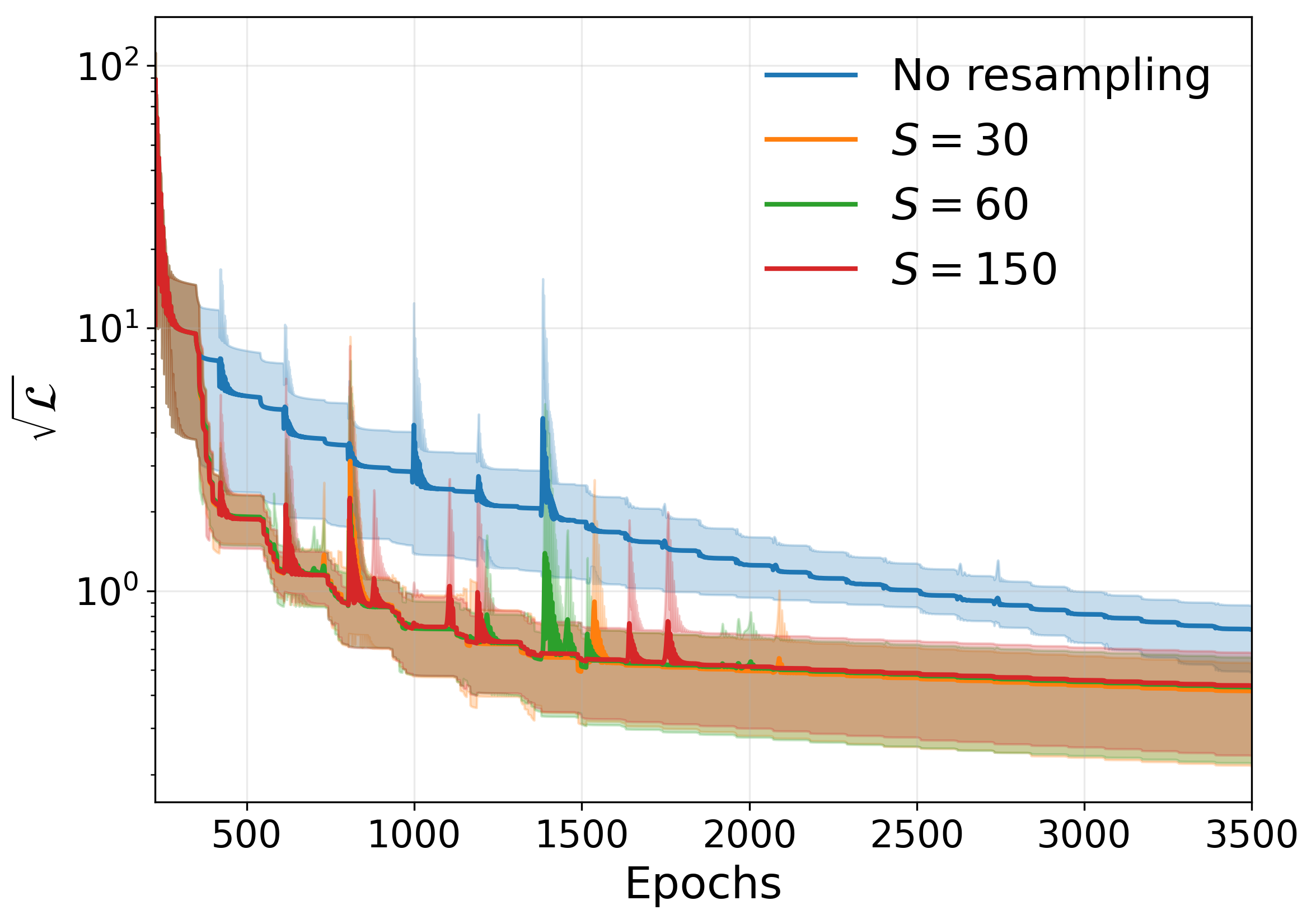}
\caption{Loss function.}
\end{subfigure}
\begin{subfigure}{0.32\linewidth}
\centering
\includegraphics[width=0.9\linewidth]{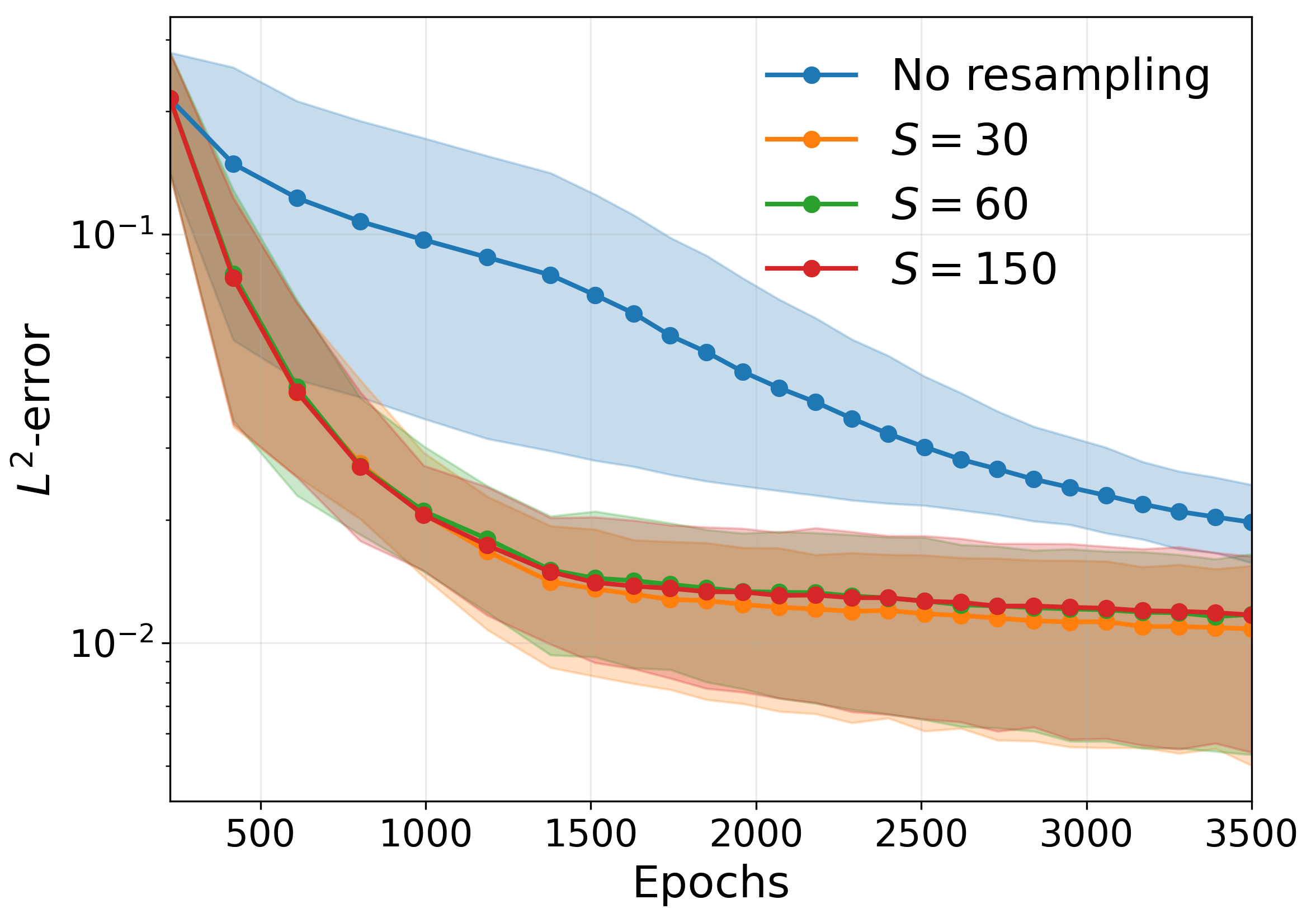}
\caption{$L^2$ error.}
\end{subfigure}
\begin{subfigure}{0.32\linewidth}
\centering
\includegraphics[width=0.9\linewidth]{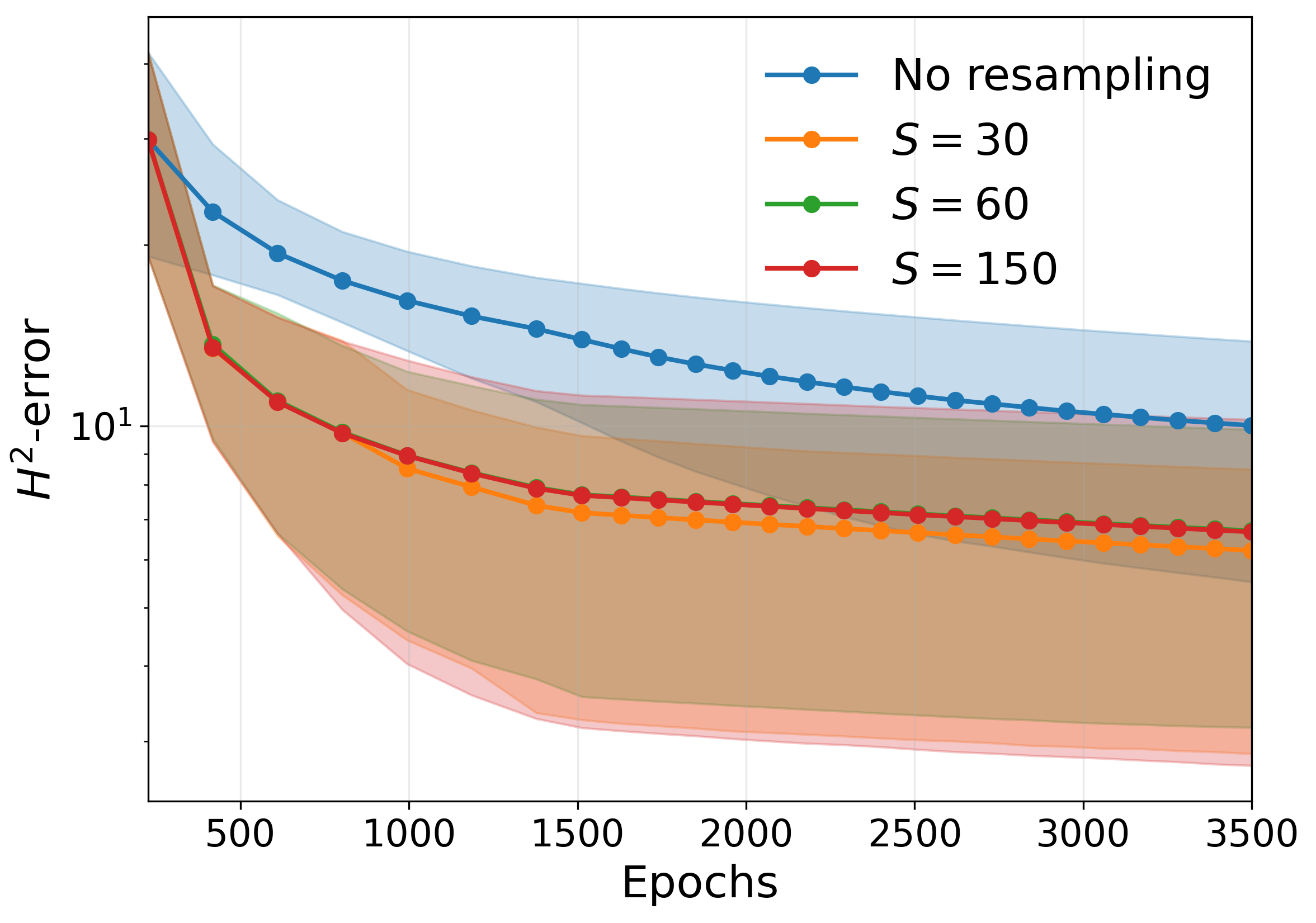}
\caption{$H^2$ error}
\end{subfigure}
\caption{First exponential test case, with $\alpha =4$. Deep Ritz method. Loss function (left), $L^2$ error (center) and $H^2$ error (right) for various values of the seed $S$.}
\label{fig:exp4_adapt}
\end{figure}

\begin{figure}[ht!]
\centering
\begin{subfigure}{0.95\linewidth}
\centering
\includegraphics[width=0.95\linewidth]{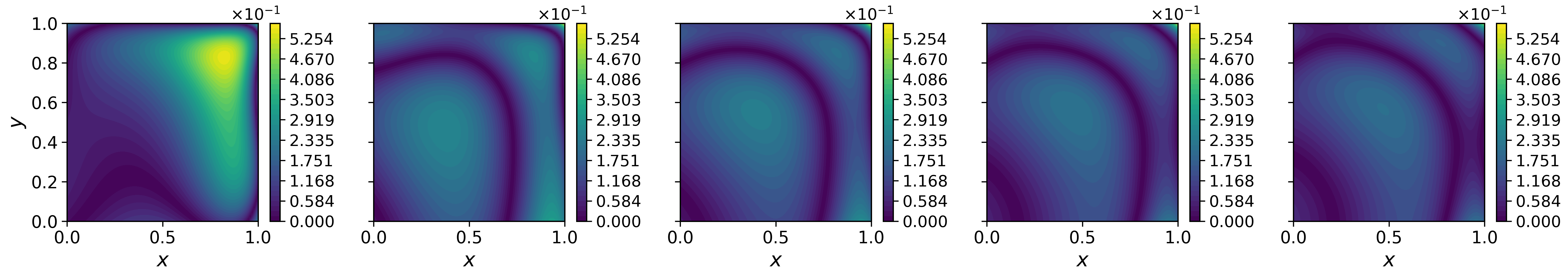}
\caption{No re-sampling.}
\end{subfigure}
\vspace{0.5cm}
\begin{subfigure}{0.95\linewidth}
\centering
\includegraphics[width=0.95\linewidth]{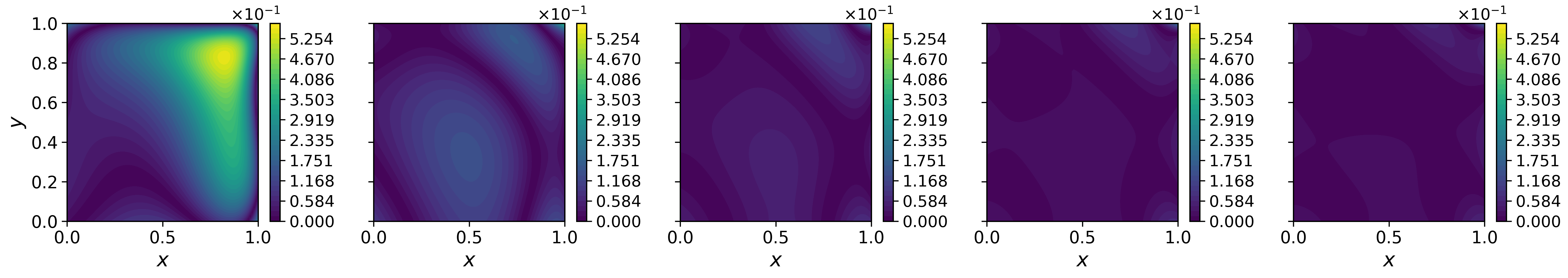}
\caption{Adaptive sampling.}
\end{subfigure}
\caption{First exponential test case, with $\alpha =4$, 3000 collocation points and 300 boundary points. 
Pointwise absolute error at the end of each splitting iteration.
Top row: no adaptive sampling; Bottom row: with adaptive sampling. }
\label{fig:exp4_pointwise}
\end{figure}

\begin{figure}[ht!]
\centering
\includegraphics[width=0.9\linewidth]{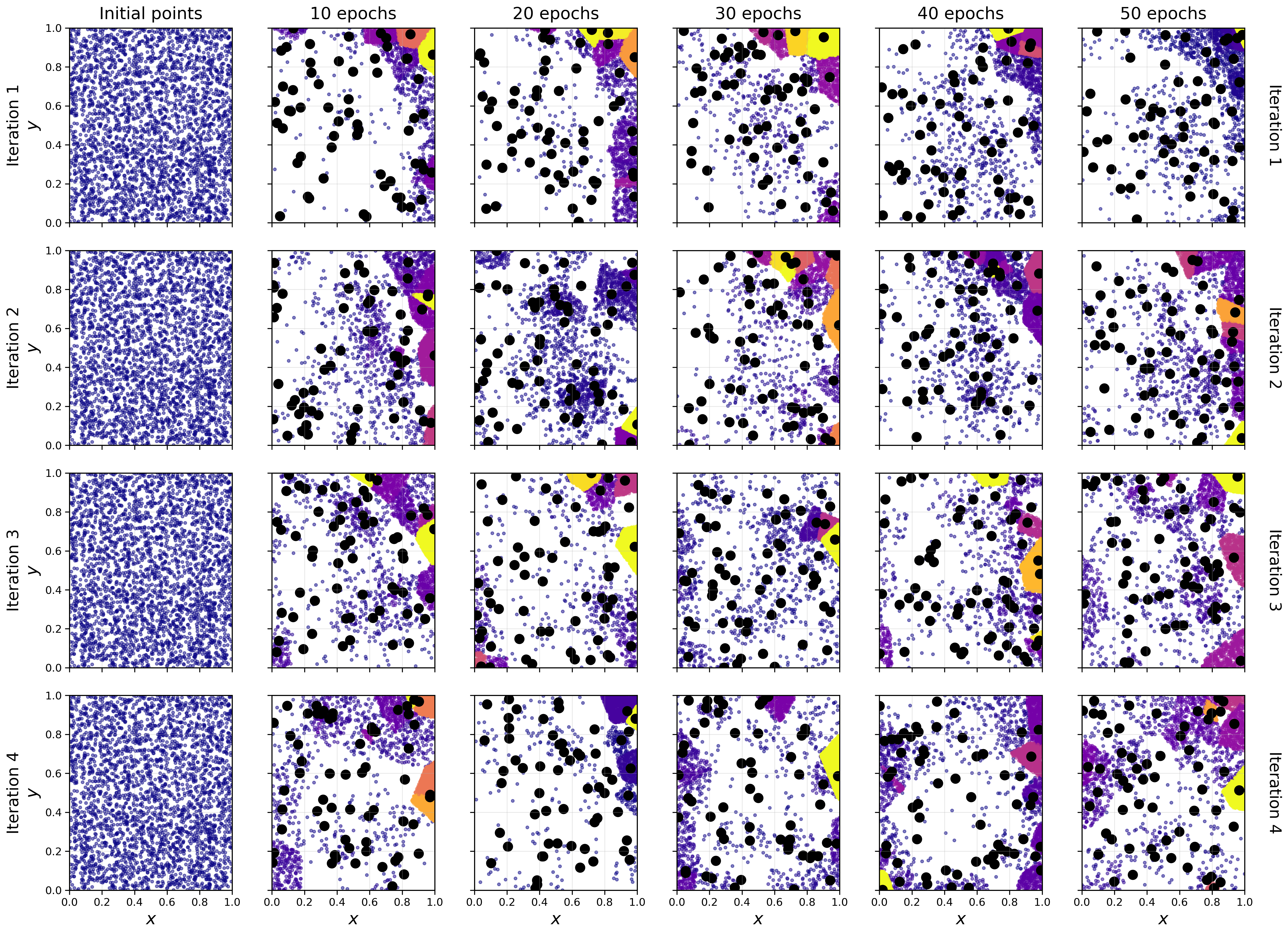}
\caption{First exponential test case, with $\alpha =4$, 3000 collocation points and 300 boundary points. 
Distribution of the resampled points (the seeds are illustrated in black).}
\label{fig:exp4_resampled}
\end{figure}

We next consider $\alpha=4$, for which the solution exhibits sharper gradients. In this more challenging regime, the standard PINNs becomes considerably less robust, whereas the Deep Ritz method retains better stability (\Cref{fig:exp_4_comp}). The benefit of adaptive sampling is more pronounced in this case. \Cref{fig:exp4_adapt,fig:exp4_pointwise} show an improvement in convergence and pointwise accuracy. The redistribution shown in \Cref{fig:exp4_resampled} is also consistent with the structure of the exact solution: the collocation points become increasingly concentrated near $(1,1)$, where the solution and its derivatives vary most rapidly.

The computational comparison in \Cref{tab:walltime-exp} confirms this behavior. For $\alpha=1$, all methods reach a relative Hessian error below $5\%$, with PINNs being the fastest. For $\alpha=4$, however, the threshold is reached in only $2/10$ PINN runs and $6/10$ standard Deep Ritz runs, compared with $9$-$10/10$ runs for the adaptive variants. Adaptive sampling also reduces the mean wall-clock time by a factor of approximately $2.2$-$2.6$ relative to standard Deep Ritz.

\newpage
\clearpage

\begin{table}[ht!]
\centering
\caption{First exponential test case: computational cost.
We report the mean wall-clock time over $10$ independent runs to reach a relative Hessian error below $5\%$; if the threshold is not reached, the time at the cap is used. The speed-up is computed relative to the standard Deep Ritz for the same value of $\alpha$.}
\label{tab:walltime-exp}
\begin{tabular}{lccc}
\toprule
Method
& Success rate
& Wall time [min]
& Speed-up \\
\midrule

\multicolumn{4}{c}{$\alpha=1$} \\
\midrule
Deep Ritz
& $10/10$
& $6.01$
& $1.00\times$ \\

PINNs
& $10/10$
& $\mathbf{1.61}$
& $\mathbf{3.72\times}$ \\

Adaptive Deep Ritz, $S=30$
& $10/10$
& $2.22$
& $2.71\times$ \\

Adaptive Deep Ritz, $S=60$
& $10/10$
& $\mathbf{2.11}$
& $\mathbf{2.84\times}$ \\

Adaptive Deep Ritz, $S=150$
& $10/10$
& $2.32$
& $2.59\times$ \\

\midrule
\multicolumn{4}{c}{$\alpha=4$} \\
\midrule
Deep Ritz
& $6/10$
& $9.21$
& $1.00\times$ \\

PINNs
& $2/10$
& $9.39$
& $0.98\times$ \\

Adaptive Deep Ritz, $S=30$
& $\mathbf{10/10}$
& $\mathbf{3.58}$
& $\mathbf{2.57\times}$ \\

Adaptive Deep Ritz, $S=60$
& $\mathbf{10/10}$
& $3.94$
& $2.33\times$ \\

Adaptive Deep Ritz, $S=150$
& $9/10$
& $4.12$
& $2.24\times$ \\
\bottomrule
\end{tabular}
\end{table}

\subsection{Monge-Amp\`ere equation in 2D : a test with sharp gradients}
We next consider the exact solution
$$u(x,y)=-\sqrt{R^2-(x^2+y^2)},\qquad R\geq\sqrt2,$$
on $\Omega=[0,1]^2$, with $f$ and $\phi$ chosen accordingly. For $R>\sqrt2$, the solution belongs to $C^\infty(\overline\Omega)$, whereas at the limiting value $R=\sqrt2$ the gradient becomes singular at $(1,1)$ and $u\notin H^2(\Omega)$.

\begin{figure}[ht!]
\centering
\begin{subfigure}{0.40\linewidth}
\centering
\includegraphics[width=0.8\linewidth]{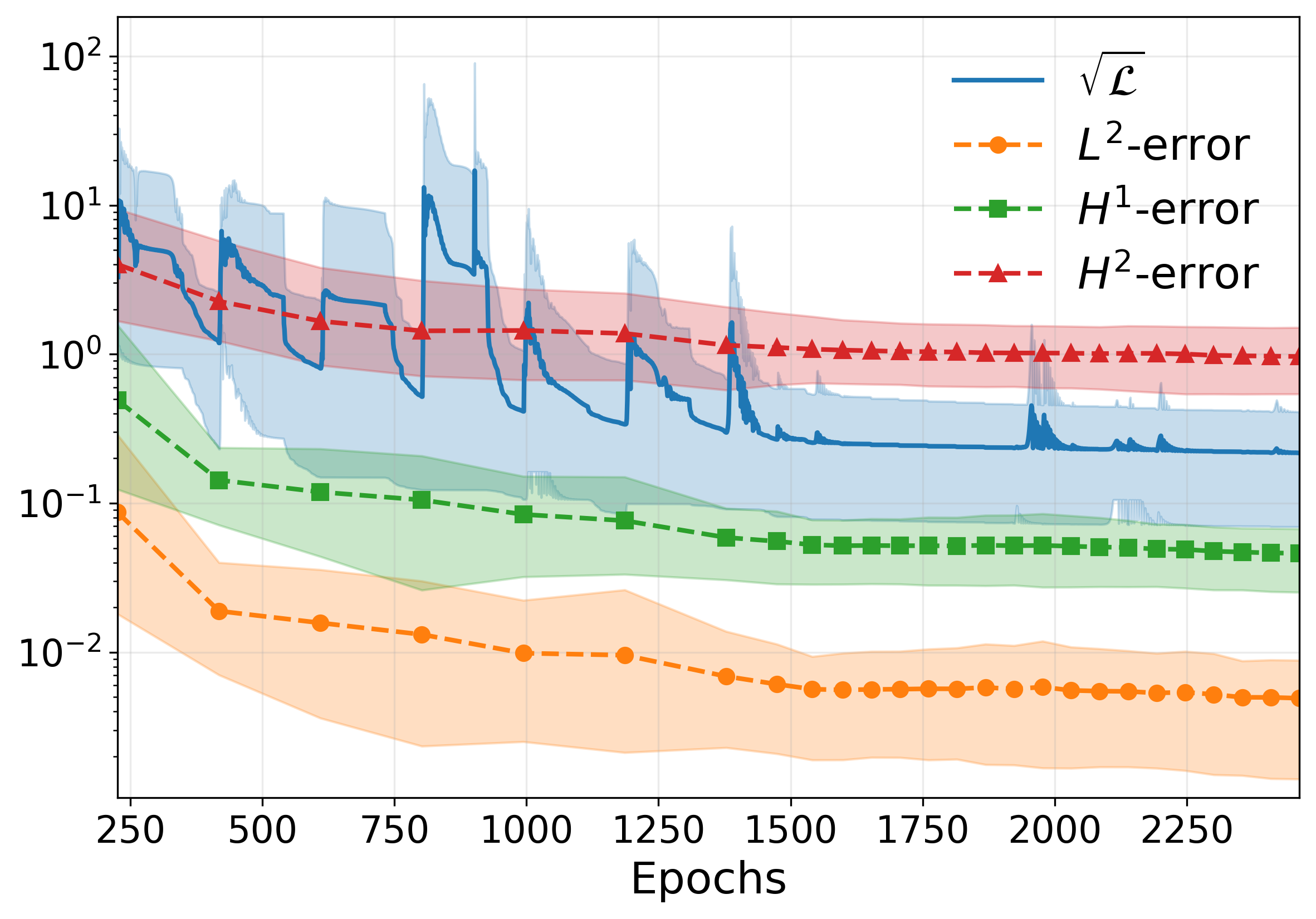}
\caption{PINNs. 3000 points.}
\end{subfigure}
\begin{subfigure}{0.40\linewidth}
\centering
\includegraphics[width=0.8\linewidth]{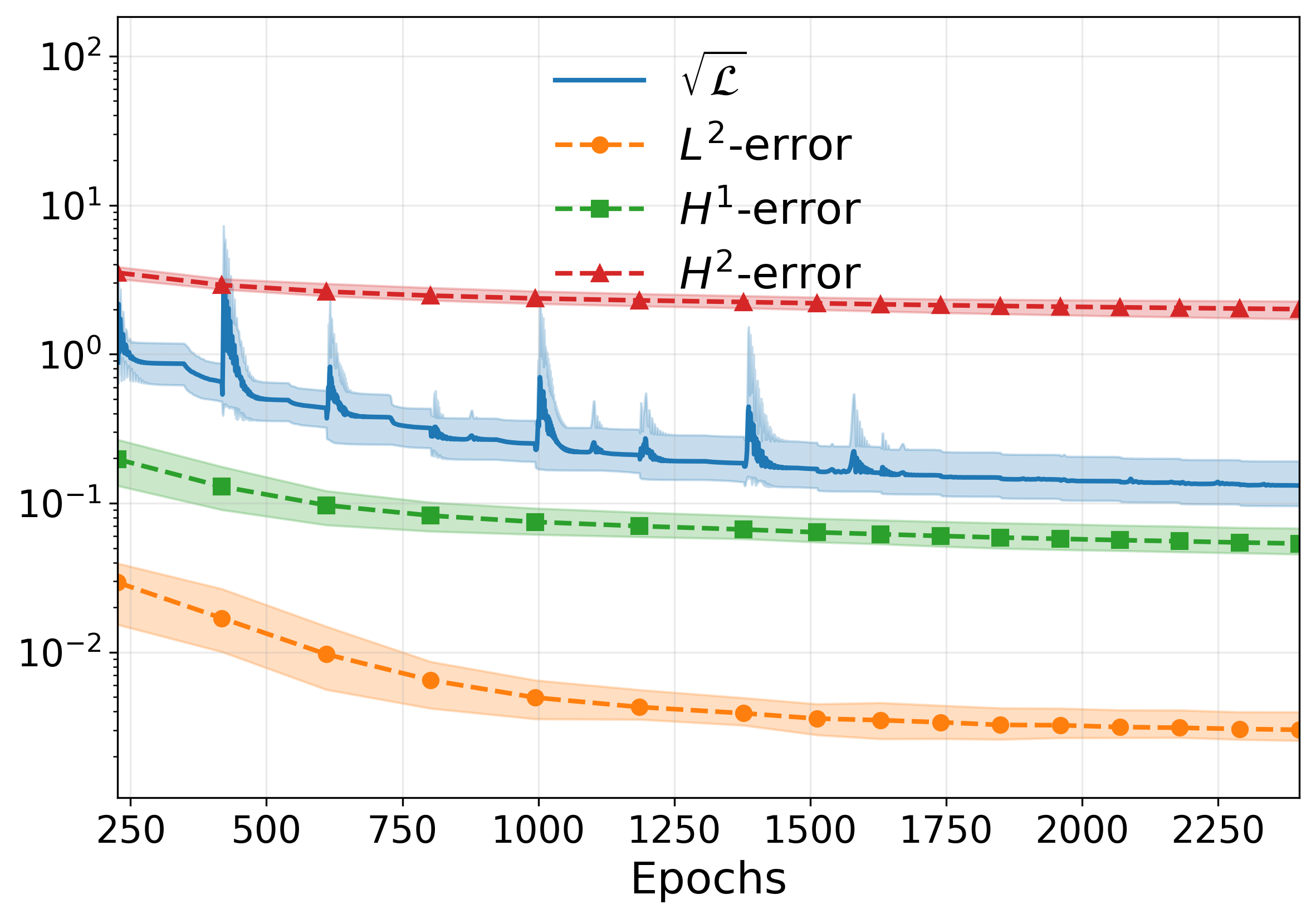}
\caption{DR. 3000 points.}
\end{subfigure}
\caption{Second test case, 3000 collocation points and 300 boundary points. 
We illustrate the loss function for each epoch and the errors between $u$ and $u_{NN}$ computed only at the end of each splitting iteration. }
\label{fig:sqrt_comp}
\end{figure}

\begin{figure}[ht!]
\centering
\begin{subfigure}{0.32\linewidth}
\centering
\includegraphics[width=0.9\linewidth]{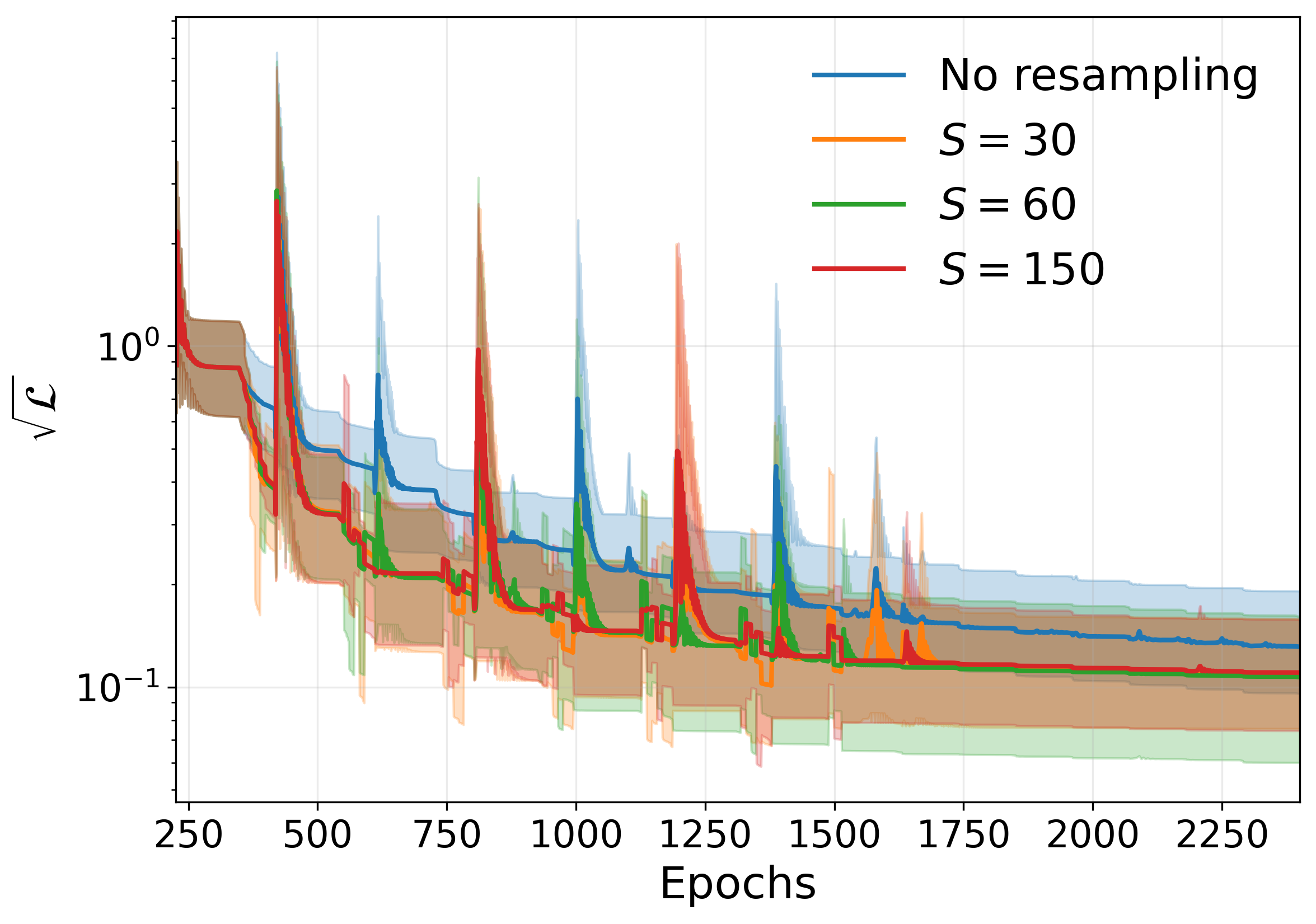}
\caption{Loss function.}
\end{subfigure}
\begin{subfigure}{0.32\linewidth}
\centering
\includegraphics[width=0.9\linewidth]{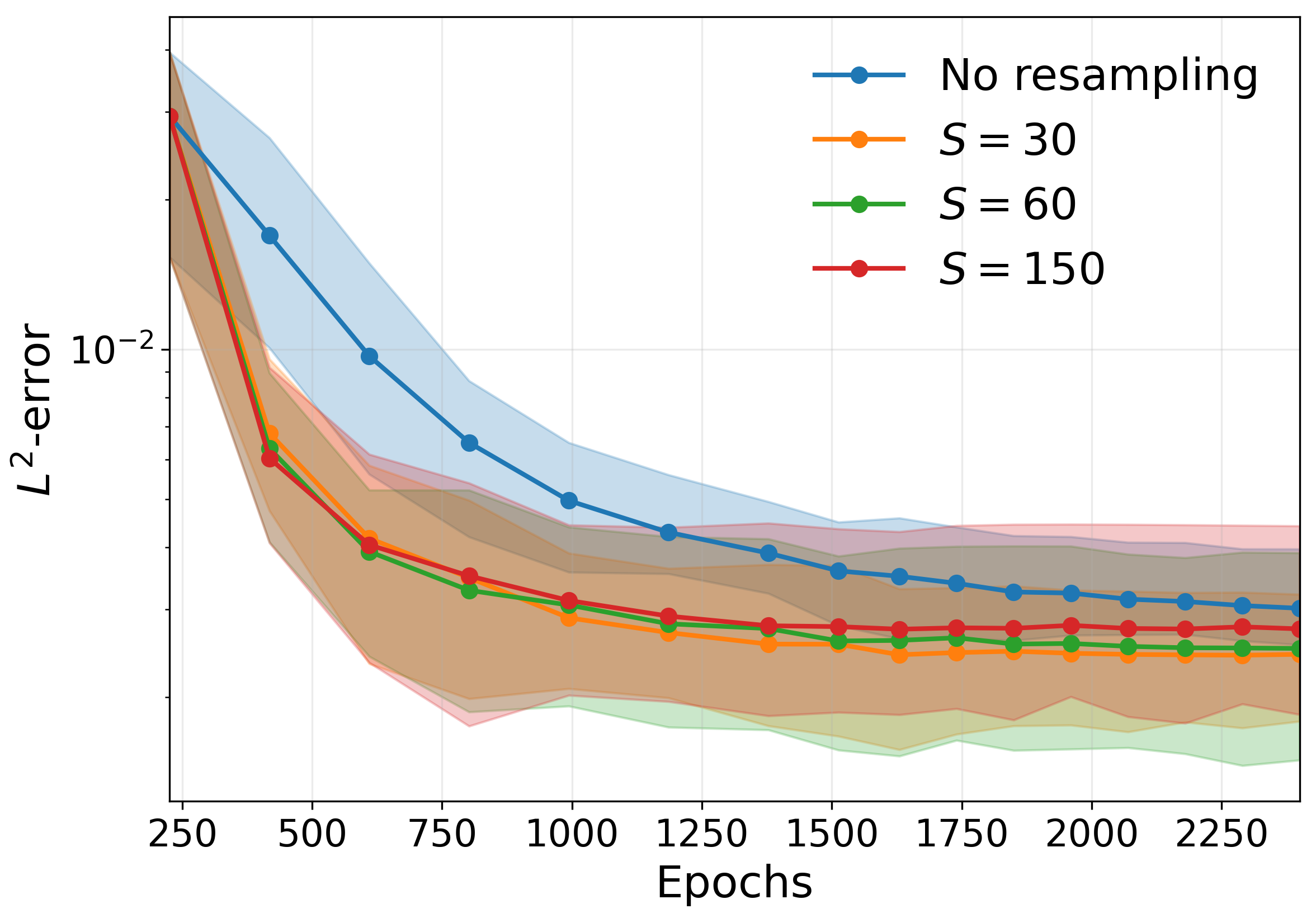}
\caption{$L^2$ error.}
\end{subfigure}
\begin{subfigure}{0.32\linewidth}
\centering
\includegraphics[width=0.9\linewidth]{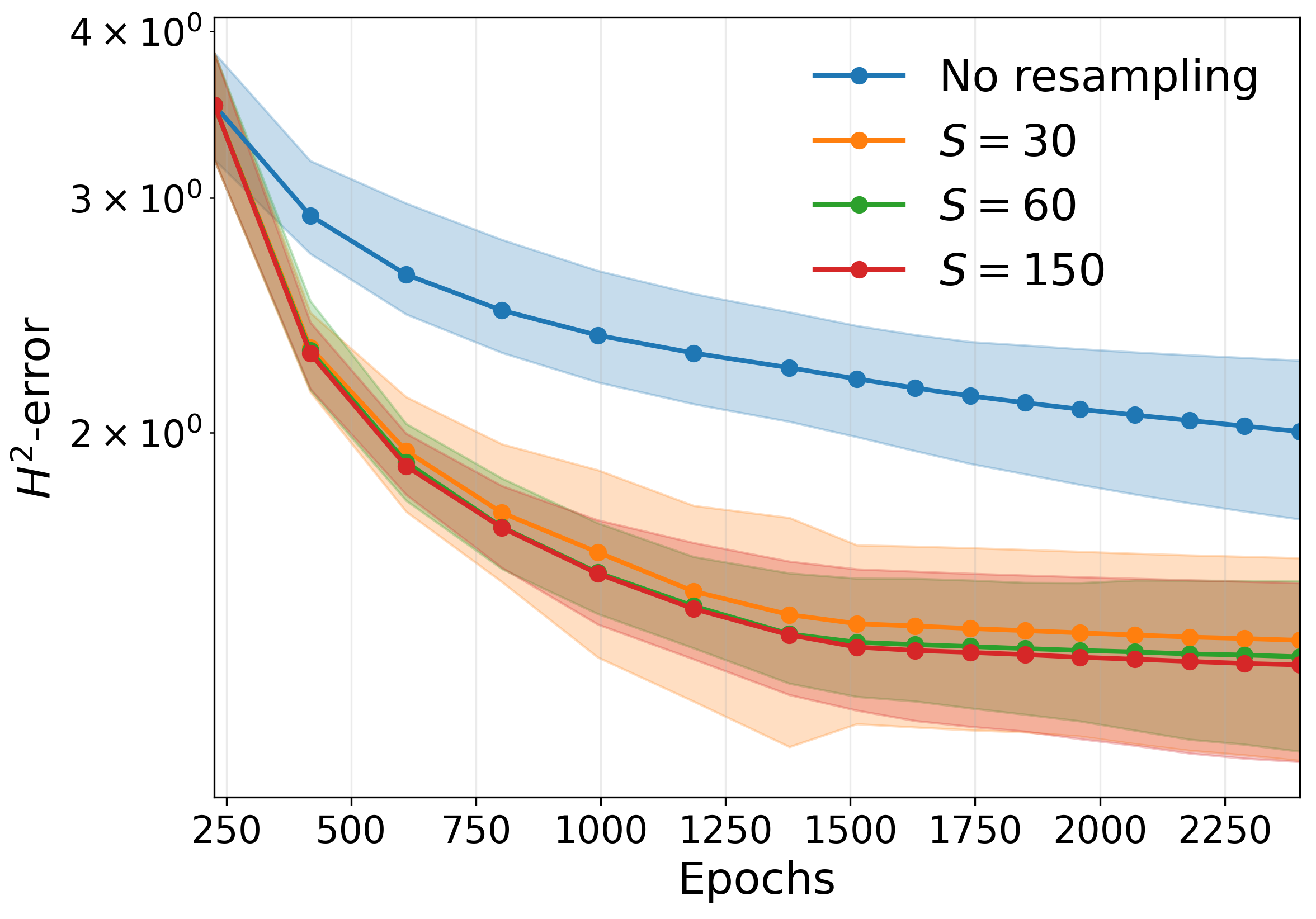}
\caption{$H^2$ error}
\end{subfigure}
\caption{Second test case. Deep Ritz method. Loss function (left), $L^2$ error (center) and $H^2$ error (right) for various values of the seed $S$. }
\label{fig:sqrt_adapt}
\end{figure}

\begin{figure}[ht!]
\centering
\begin{subfigure}{0.95\linewidth}
\centering
\includegraphics[width=0.95\linewidth]{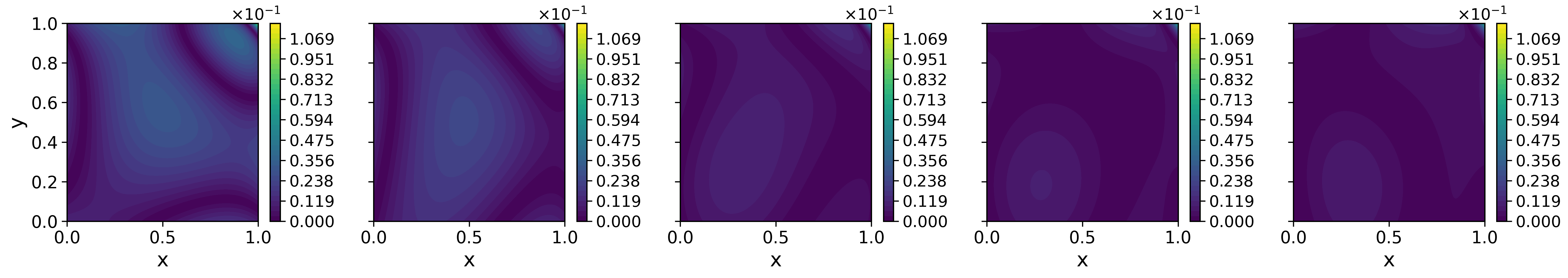}
\caption{Without adaptive sampling.}
\end{subfigure}
\vspace{0.5cm}
\begin{subfigure}{0.95\linewidth}
\centering
\includegraphics[width=0.95\linewidth]{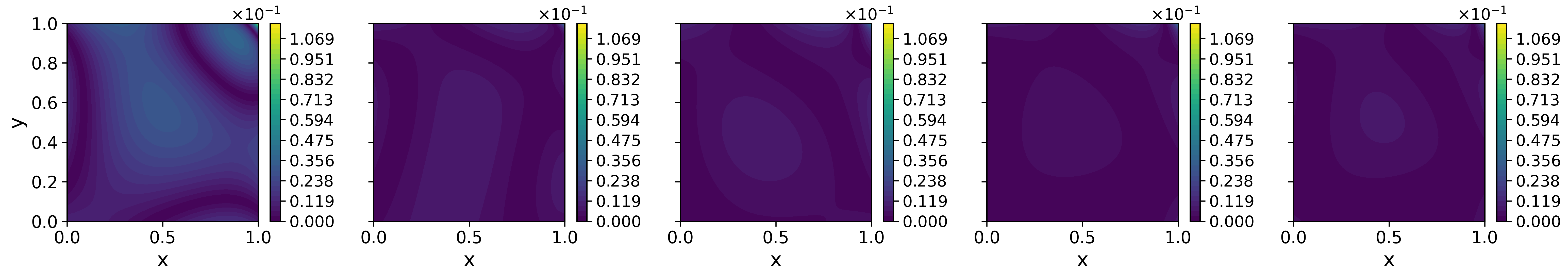}
\caption{With adaptive sampling.}
\end{subfigure}
\caption{
Second test case, $3000$ collocation points and $300$ boundary points. Deep Ritz method. 
Pointwise absolute error at the end of each splitting iteration. 
Top row: without adaptive sampling; bottom row: with adaptive sampling.
}
\label{fig:sqrt_pointwise}
\end{figure}

We take $R=\sqrt{2}+0.01$. As shown in \Cref{fig:sqrt_comp}, PINNs and Deep Ritz give comparable accuracy while Deep Ritz shows less variability. Adaptive sampling further improves the Deep Ritz approximation, with little sensitivity to the number of seeds (\Cref{fig:sqrt_adapt}). The pointwise errors in \Cref{fig:sqrt_pointwise} show that adaptive sampling both reduces the error and distributes it more uniformly over the domain.

\subsection{Monge-Amp\`ere equation in 2D : a test on the unit disk}
\begin{figure}[ht!]
\centering
\includegraphics[width=0.4\linewidth]{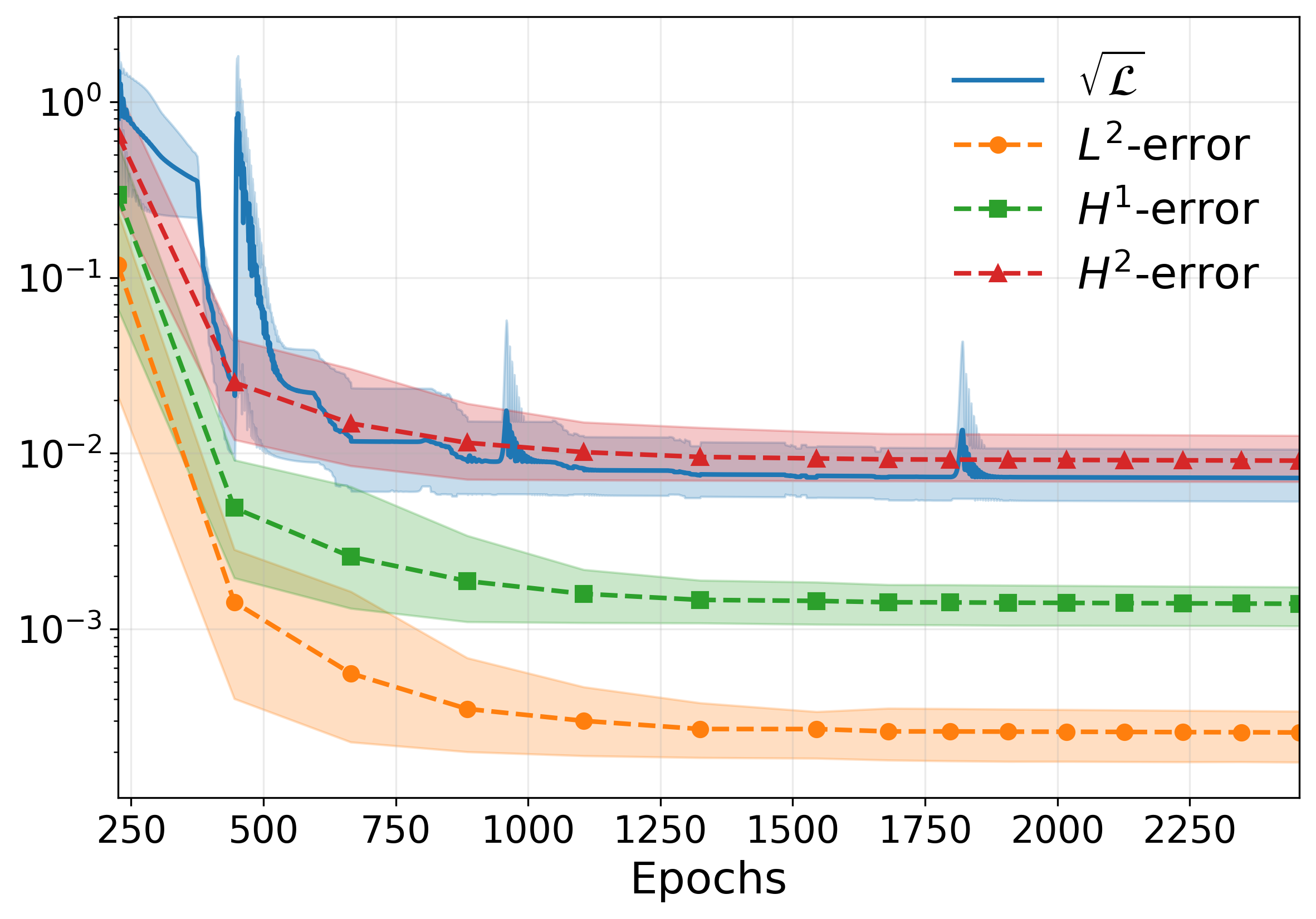}
\caption{A test case on the unit disk, with $3000$ collocation points and 300 boundary points, together with adaptive sampling. 
We illustrate the loss function for each epoch and the errors between $u$ and $u_{NN}$ computed only at the end of each splitting iteration. 
}
\label{fig:gauss_disk}
\end{figure}

\begin{figure}[ht!]
\centering
\begin{subfigure}{0.95\linewidth}
\centering
\includegraphics[width=0.95\linewidth]{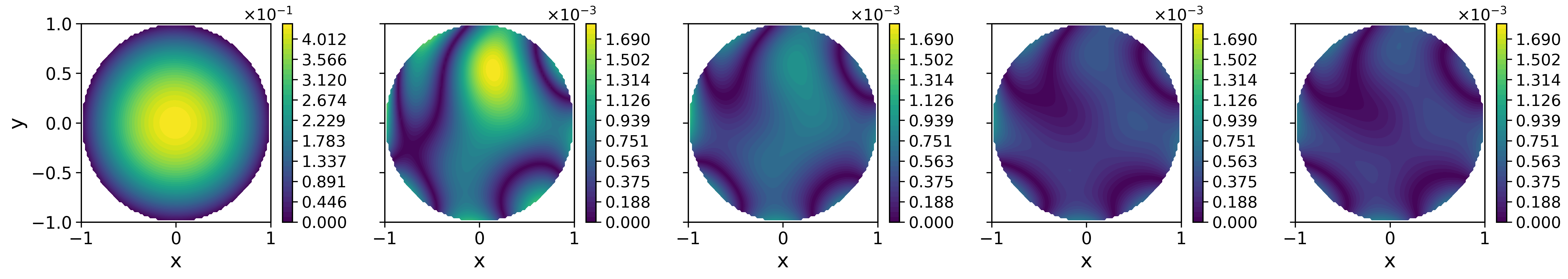}
\caption{Without adaptive sampling.}
\end{subfigure}
\vspace{0.5cm}
\begin{subfigure}{0.95\linewidth}
\centering
\includegraphics[width=0.95\linewidth]{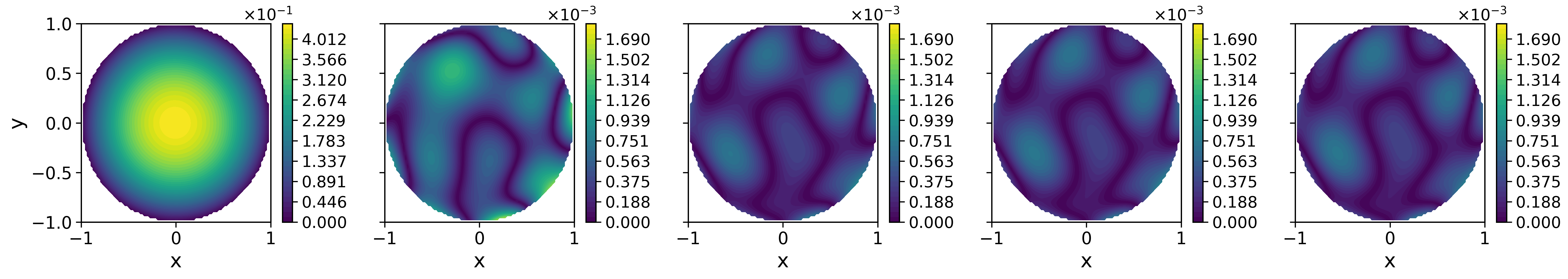}
\caption{With adaptive sampling.}
\end{subfigure}
\caption{A test case on the unit disk, with $3000$ collocation points and $300$ boundary points. Deep Ritz method. 
Pointwise absolute error at the end of each splitting iteration. 
Top row: without adaptive sampling; bottom row: with adaptive sampling.
}
\label{fig:gauss_disk_pointwise}
\end{figure}

In order to highlight that the Deep Ritz method can handle curved boundaries and nearly elliptic cases, we consider the case 
$$
\Omega= \left\{ (x,y) \in \mathbb{R}^2 : x^2 + y^2 < 1 \right\}
$$ 
with $f = 1$ and $g = 0$. This test case admits the exact solution 
$$
u(x, y) = \frac12 \left( x^2 + y^2  -1 \right)
$$
\Cref{fig:gauss_disk} shows the convergence of the Deep Ritz method, while \Cref{fig:gauss_disk_pointwise} shows that adaptive sampling reduces and more uniformly distributes the pointwise error. This example also illustrates that the mesh-free formulation accommodates curved boundaries without any modification of the method.

\subsection{Extension to the Pucci's equation in 2D}

For the Pucci equation \eqref{intro_pucci}, the local nonlinear projection is replaced by the solver of \cite{pucci}, and, since convexity is no longer required, the ICNN is replaced by a standard feedforward network with softplus activation. We set $\lambda=1000$.

\begin{figure}[ht!]
\centering
\begin{subfigure}{0.32\linewidth}
\centering
\includegraphics[width=0.9\linewidth]{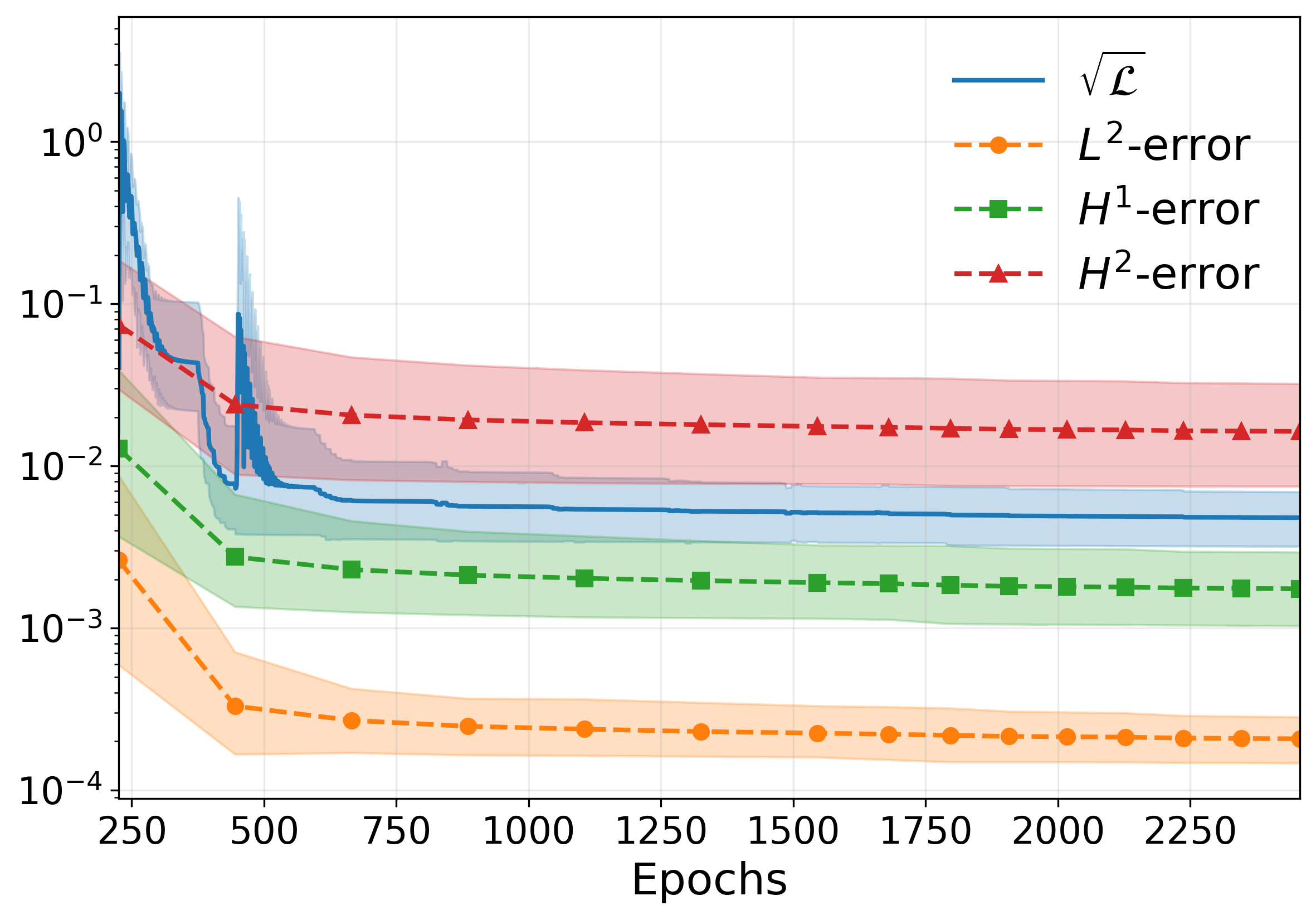}
\caption{$\alpha=2$.}
\end{subfigure}
\begin{subfigure}{0.32\linewidth}
\centering
\includegraphics[width=0.9\linewidth]{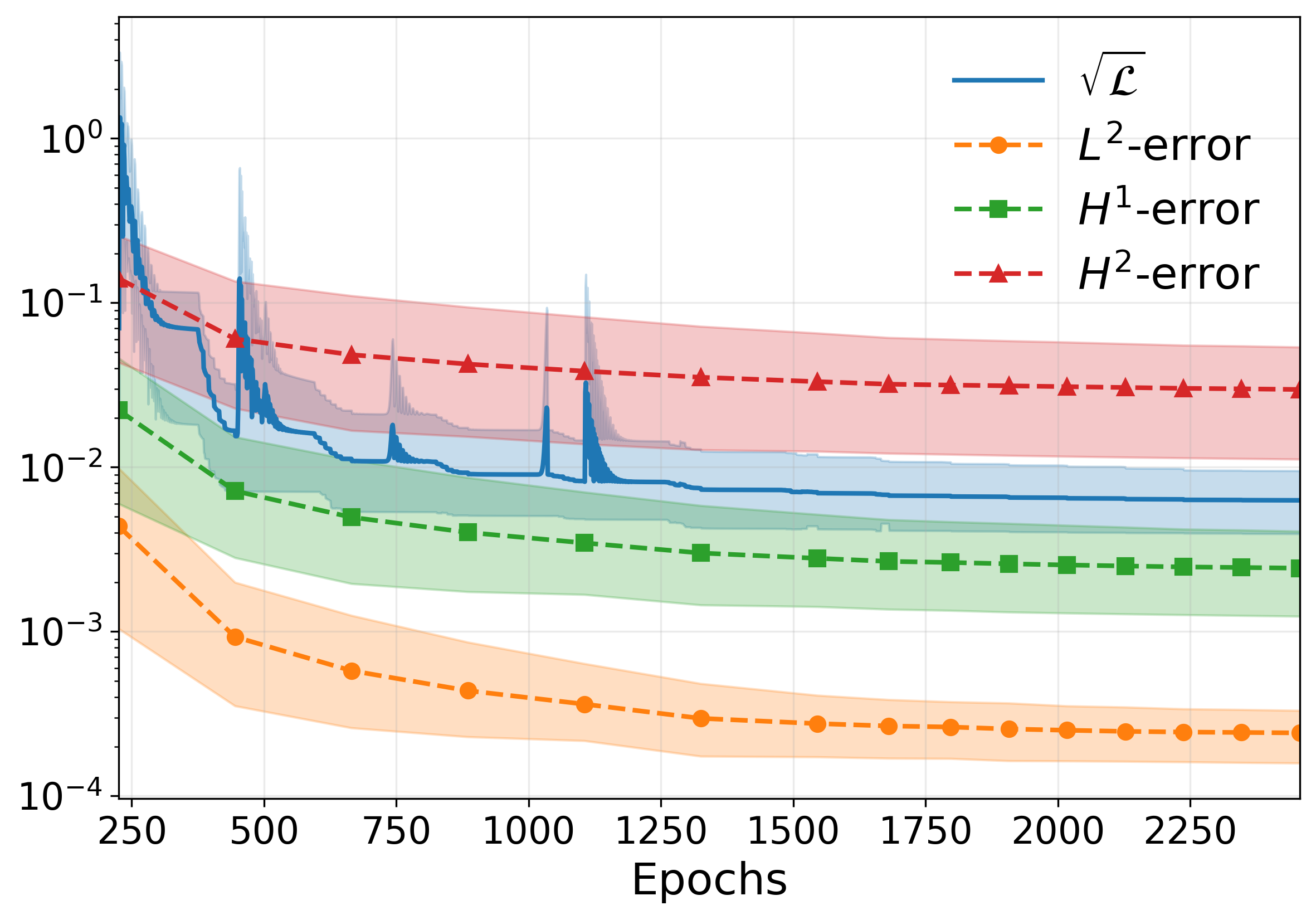}
\caption{$\alpha=3$.}
\end{subfigure}
\begin{subfigure}{0.32\linewidth}
\centering
\includegraphics[width=0.9\linewidth]{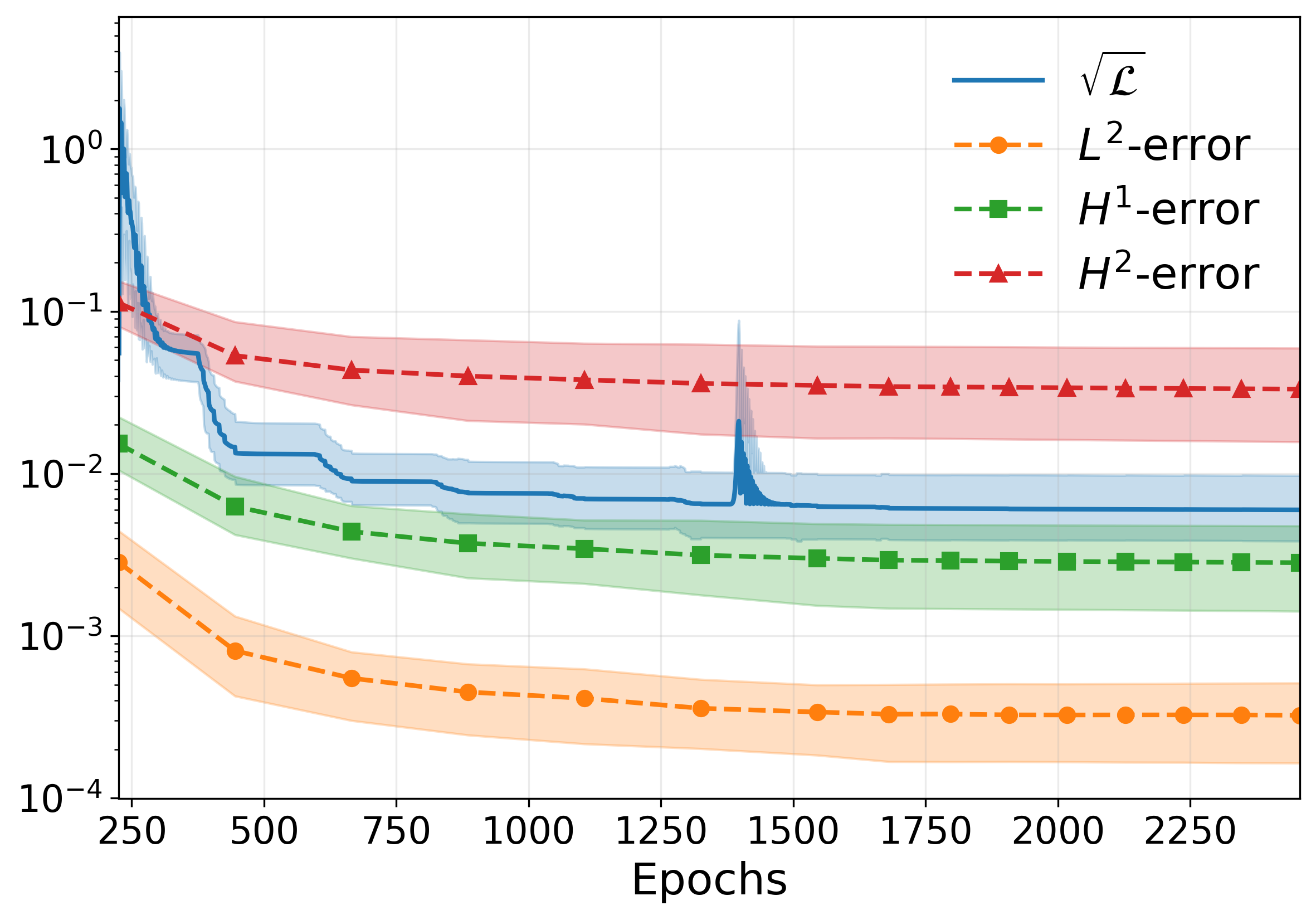}
\caption{$\alpha=5$.}
\end{subfigure}
\caption{Test case for the Pucci's equation. Convergence results when using the Deep Ritz method with $3000$ collocation points and $300$ boundary points.}
\label{fig:pucci}
\end{figure}

On $\Omega=[0,1]^2$, we consider the exact solution
$$u(x,y)=-((x+1)^2+(y+1)^2)^{(1-\alpha)/2}$$
with the corresponding data $f$ and $\phi$. The solution is smooth on $\overline\Omega$ for every $\alpha>1$, while increasing $\alpha$ produces progressively sharper variations. As shown in \Cref{fig:pucci}, the method remains stable and accurate for $\alpha=2,3,5$, illustrating that the proposed computational framework extends beyond the Monge-Amp\`ere equation.

\subsection{Extension to the Gauss curvature equation in 2D: The Minkowski problem}

\begin{figure}[ht!]
\centering
\includegraphics[width=0.4\linewidth]{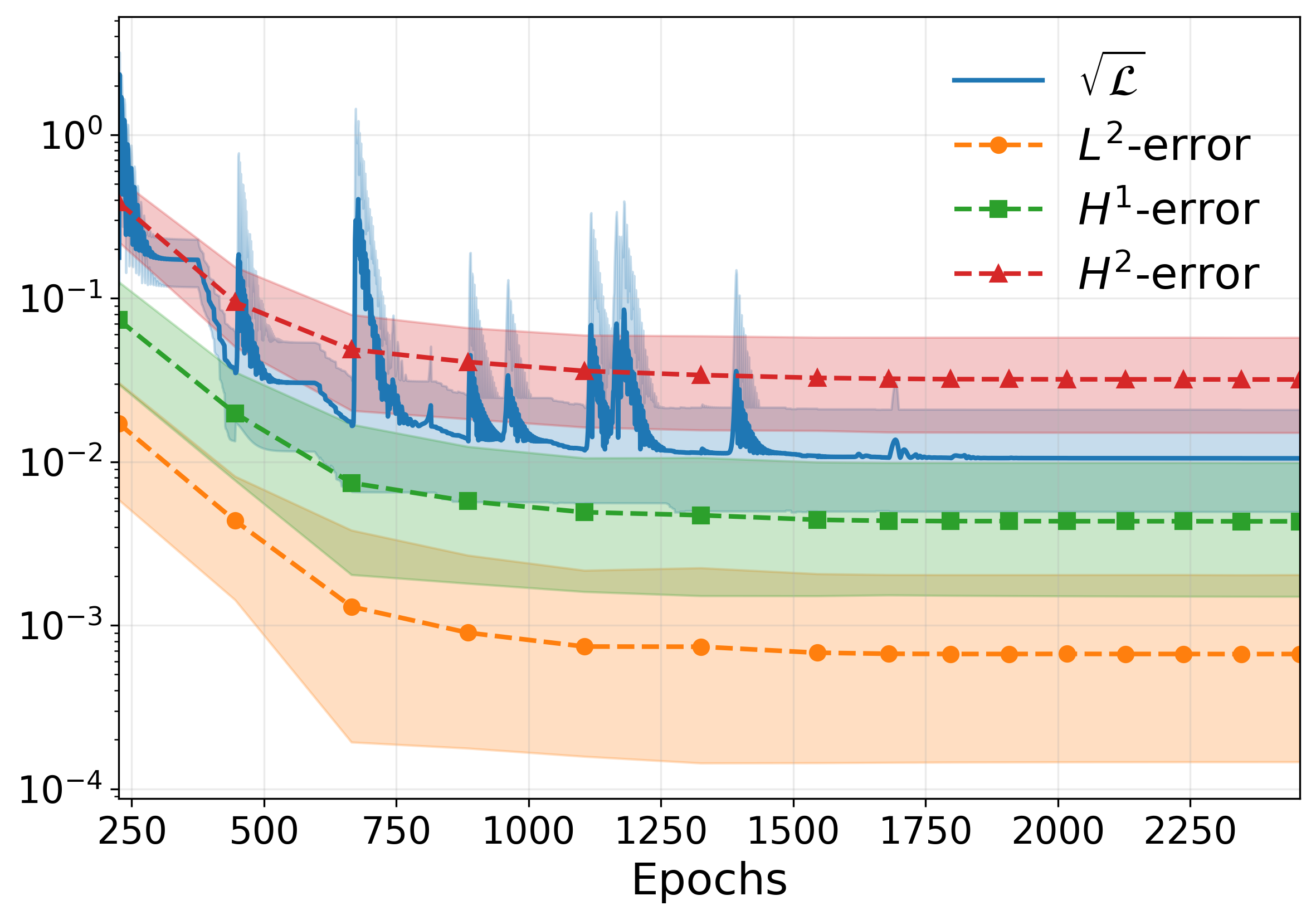}
\caption{
Gauss curvature equation in 2D. Convergence results with the adaptive sampling procedure, $3000$ collocation points and 300 boundary points.}
\label{fig:gauss_curv}
\end{figure}

\begin{figure}[ht!]
\centering
\begin{subfigure}{0.95\linewidth}
\centering
\includegraphics[width=0.95\linewidth]{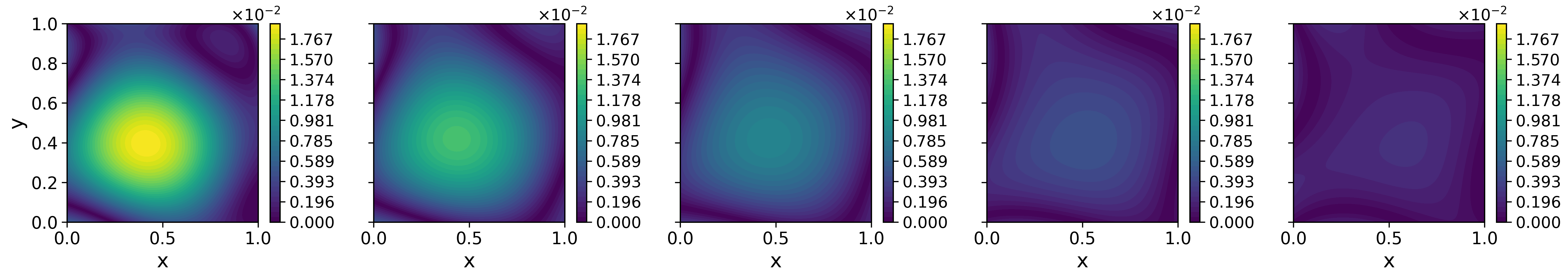}
\caption{Without adaptive sampling.}
\end{subfigure}
\vspace{0.5cm}
\begin{subfigure}{0.95\linewidth}
\centering
\includegraphics[width=0.95\linewidth]{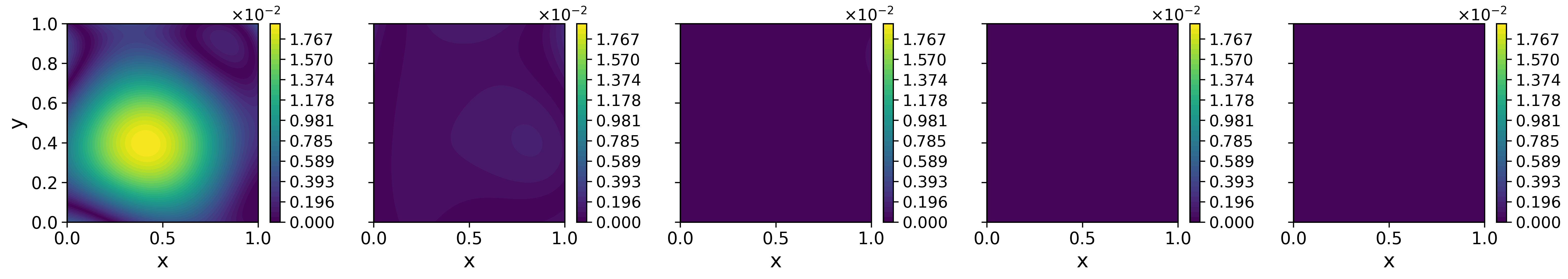}
\caption{With adaptive sampling.}
\end{subfigure}
\caption{
Gauss curvature equation in 2D, with $3000$ collocation points and $300$ boundary points. Deep Ritz method. 
Pointwise absolute error at the end of each splitting iteration.
Top row: without adaptive sampling; bottom row: with adaptive sampling.
}
\label{fig:gauss_curv_pointwise}
\end{figure}

We finally consider the Gauss curvature equation on $\Omega=[0,1]^2$, find $u$ satisfying
\begin{equation*}
\begin{cases}
\det D^2u=K(x)(1+|\nabla u|^2)^2,&\text{in }\Omega,\\
u=g,&\text{on }\partial\Omega,
\end{cases}
\end{equation*}
where, with $\mathbf b=(1/2,1/2)$,
$$g(\mathbf x)=|\mathbf x-\mathbf b|^2,\qquad K(\mathbf x)= \frac{4}{(1+4|\mathbf x-\mathbf b|^2)^2}.$$
The exact solution is $u(\mathbf x)=|\mathbf x-\mathbf b|^2.$ Despite the additional gradient dependence in the nonlinear term, the convergence behavior remains comparable to that observed for the previous equations (\Cref{fig:gauss_curv}). Moreover, \Cref{fig:gauss_curv_pointwise} shows that adaptive sampling substantially
reduces and homogenizes the pointwise error.

\section{Numerical results for optimal transport problems}\label{sec:numres_ext}

Let us turn now to the optimal transport Monge-Amp\`ere problem \eqref{eq:Jacobian} \eqref{eq:transportBC}. 
In a first step, we tackle test cases that transport uniform distributions from one domain to another.
Then we discuss cases in which the distributions are Gaussians, and mimic the physical transport of mounts of material. 

\subsection{Optimal transport Monge-Amp\`ere preliminary test case: transporting a disk domain into an ellipse}
The aim is to transport a constant probability density function supported in the unit disk domain $ \mathcal{X} \subset \mathbb{R}^2$ to another constant probability density function whose support is an ellipse domain $ \mathcal{Y} \subset \mathbb{R}^2$.
More precisely, let us thus define
$$ 
\mathcal{X}=\{x=(x_1,x_2)\in\mathbb{R}^2:\,  x_1^2+x_2^2 < 1\}, \quad \textrm{ and  } \quad 
\mathcal{Y} = \left\{y=(y_1,y_2) \in \mathbb{R}^2:\, \frac{\left(y_1 - 3.5 \right)^2}{(2)^2}+\frac{y_2^2}{(0.5)^2} < 1\right\}.
$$ 
The corresponding generalized Monge-Amp\`ere problem \eqref{eq:Jacobian} \eqref{eq:transportBC} is considered here with 
$f(x) = \frac{1}{\pi}\chi_{\mathcal{X}}(x)$, $g(y) = \frac{1}{\pi}\chi_{\mathcal{Y}}(y)$, where $\chi_{A}$ denotes the characteristic function of the domain $A$. 
With this data, the exact solution to \eqref{eq:Jacobian} \eqref{eq:transportBC} is 

$$
u_{ex}(x) = x_1^2 + \frac14 x_2^2 + \frac72 x_1 - \frac{5}{16}, 
\qquad 
\nabla u_{ex}(x) 
=\begin{pmatrix} 2x_1 + \frac{7}{2} \\ \frac{1}{2}x_2 \end{pmatrix}, 
\qquad 
x = (x_1,x_2)\in \mathcal{X}.
$$

The initial guess $u^0$ is computed in order to guarantee that  $\nabla u^0$ is the identity tensor. 
We discretize the loss function \eqref{eq:tildeEb} with $1000$ collocation points and $1000$ boundary points.

\begin{figure}[ht!]
\centering
\begin{subfigure}{0.40\linewidth}
\centering
\includegraphics[width=0.8\linewidth]{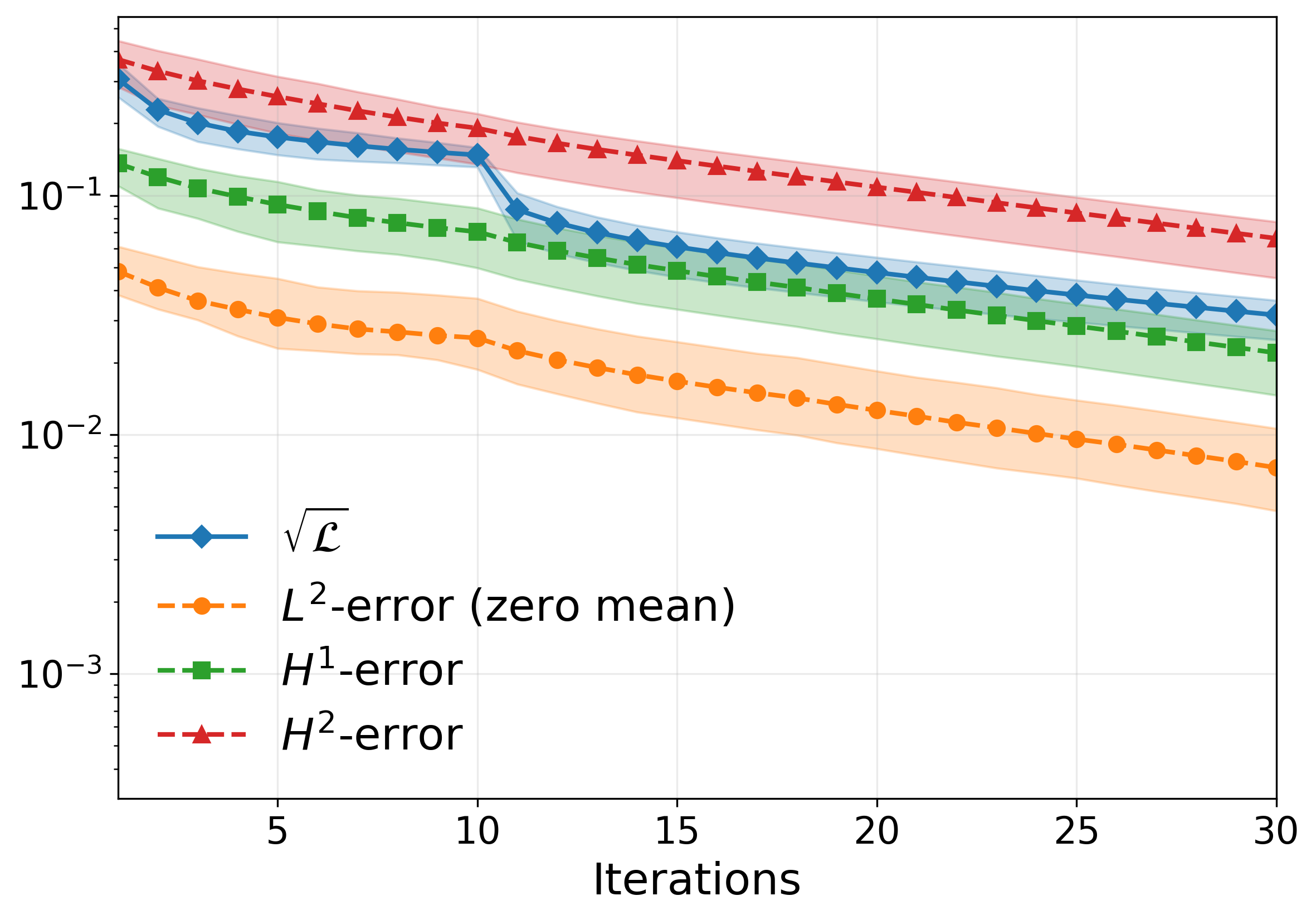}
\caption{No resampling.}
\end{subfigure}
\begin{subfigure}{0.40\linewidth}
\centering
\includegraphics[width=0.8\linewidth]{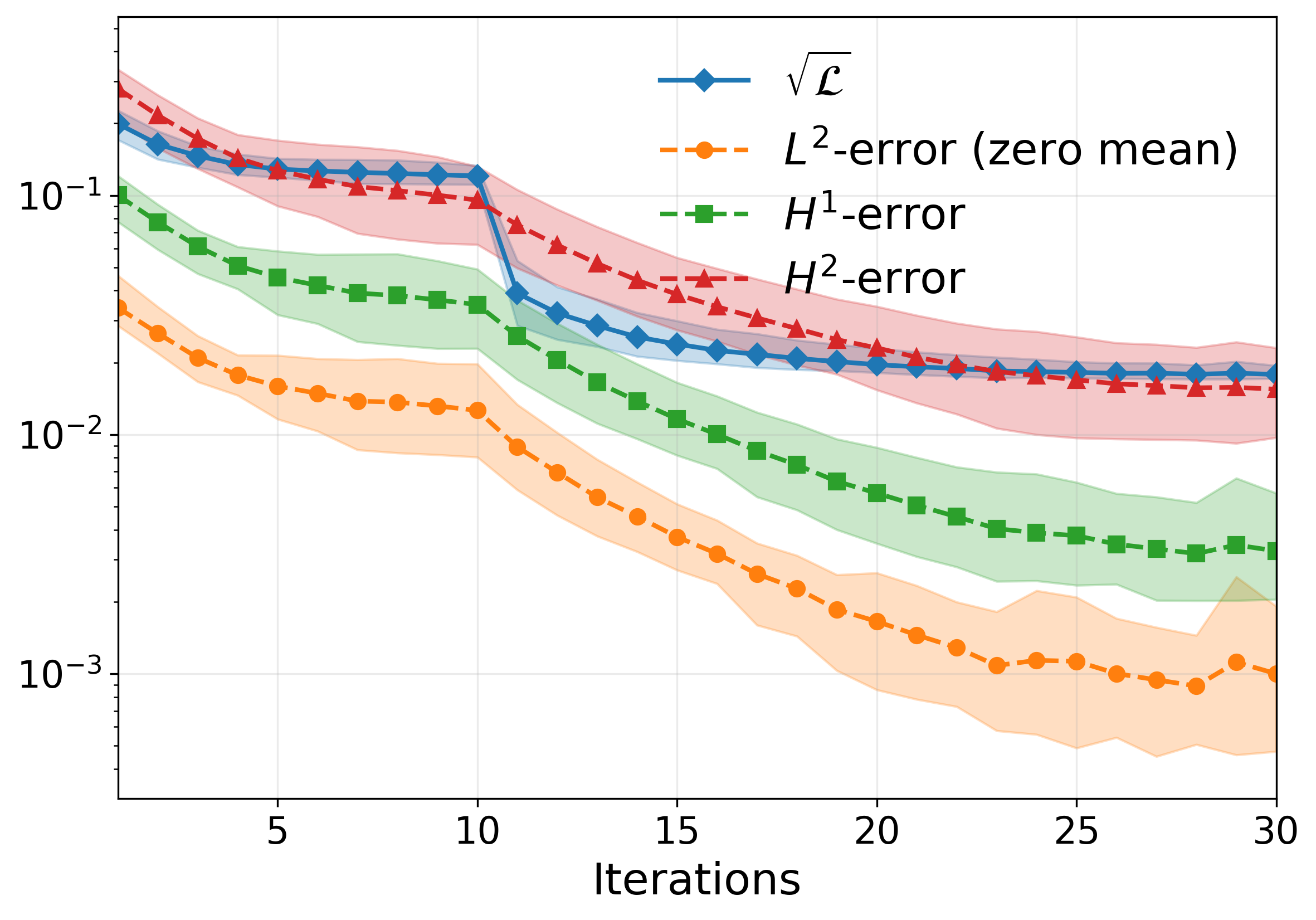}
\caption{Adaptive sampling.}
\end{subfigure}
\caption{Optimal transport Monge-Amp\`ere problem. Disk domain into an ellipse. We illustrate the loss function and the errors between $u$ and $u_{NN}$ computed at the end of each splitting iteration. 
}
\label{fig:ot_comp}
\end{figure}

\Cref{fig:ot_comp} compares the convergence of the algorithm with and without the adaptive sampling. The plots show that the adaptive sampling  strategy is faster as also observed in the other test cases. 

\begin{figure}[ht!]
\centering
\includegraphics[width=0.95\linewidth]{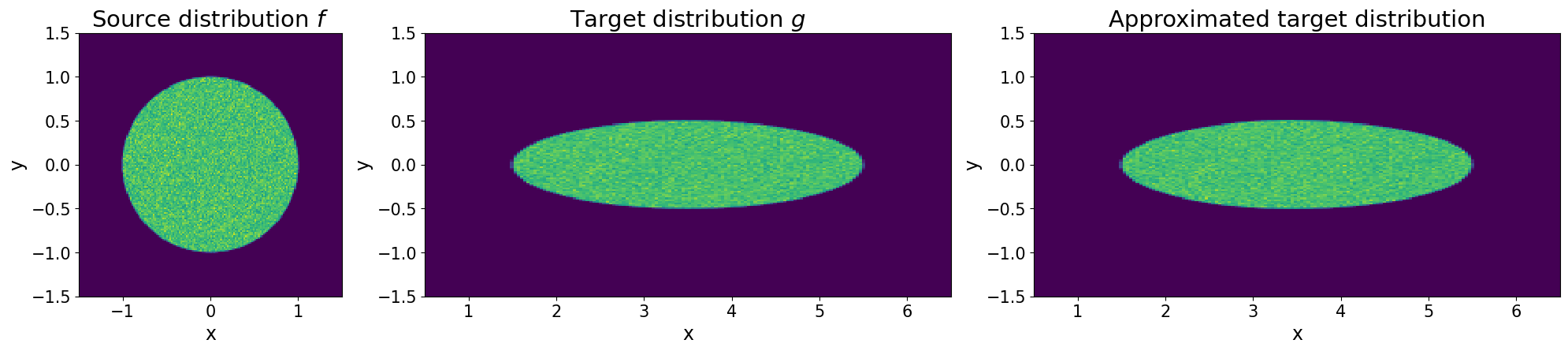}
\caption{
Optimal transport Monge-Amp\`ere problem. Disk domain into an ellipse. 
Visualization of the transport map with histograms based on $10^6$ sampling points. 
Results after $30$ iterations of the splitting algorithm. 
Left: Source distribution $f$; 
middle: Exact target distribution $g$; 
right: approximated target distribution $(\nabla u_{NN})_{\#}(f)$. 
}
\label{fig:otfirst}
\end{figure}

\Cref{fig:otfirst} visualizes the optimal transport map by showing the location of $10^6$ uniformly sampled points in the unit ball and their image into the ellipse, both with the exact solution and with the approximate map. It shows that exact and approximate maps give very similar solutions.

\subsection{Optimal transport Monge-Amp\`ere: transporting Gaussian distributions}
We consider three optimal transport problems on \(\Omega=[0,1]^2\), with increasingly complex source and target distributions. In all cases, the transport map is approximated by $\nabla u_{NN}$. To visualize the transport, we sample $10^6$ points from the source distribution and display their push-forward distribution at several splitting iterations, with and without adaptive sampling.

We first consider a Gaussian source $f_1$ transported to the uniform distribution $g_{\rm unif}$, namely
\begin{equation}\label{eq:gaussian1}
f_1(x)
=
c_1
\exp\left(
-\frac{|x-x_0|^2}{2\sigma^2}
\right)
\chi_{\Omega}(x),
\qquad
x_0=(0.25,0.75),
\qquad
\sigma^2=0.25,
\end{equation}
and $$g_{\rm unif}(x)=\chi_{\Omega}(x),$$
where $c_1$ is chosen so that $f_1$ has unit mass. 

We next consider the anisotropic bimodal source
\begin{equation}\label{eq:gaussian2}
\begin{aligned}
f_2(x)=c_2\Bigg[&\exp\left(-\frac{(x_1-0.5)^2}{2\sigma_{xx}^2}-\frac{(x_2-0.2)^2}{2\sigma_{yy}^2}\right)\\
&+\exp\left(-\frac{(x_1-0.5)^2}{2\sigma_{xx}^2}-\frac{(x_2-0.8)^2}{2\sigma_{yy}^2}\right)\Bigg]\chi_{\Omega}(x),
\end{aligned}
\end{equation}
with $\sigma_{xx}^2=0.25$, $\sigma_{yy}^2=0.015625$. Transporting $f_2$ to $g_{\rm unif}$ provides a more challenging example with two spatially separated sources. 
Finally, we transport the same bimodal source $f_2$ to the centered Gaussian target
\begin{equation}\label{eq:gaussian1bis}
g_{\rm G}(x)=c_g\exp\left(-\frac{(x_1-0.5)^2+(x_2-0.5)^2}{2\sigma_g^2}
\right)\chi_{\Omega}(x),\qquad \sigma_g=0.2,
\end{equation}
where $c_g$ is the corresponding normalization constant. This is the most demanding of the three tests and requires a larger number of splitting iterations. 

\begin{figure}[ht!]
\centering
\begin{subfigure}{0.95\linewidth}
\centering
\includegraphics[width=0.99\linewidth]{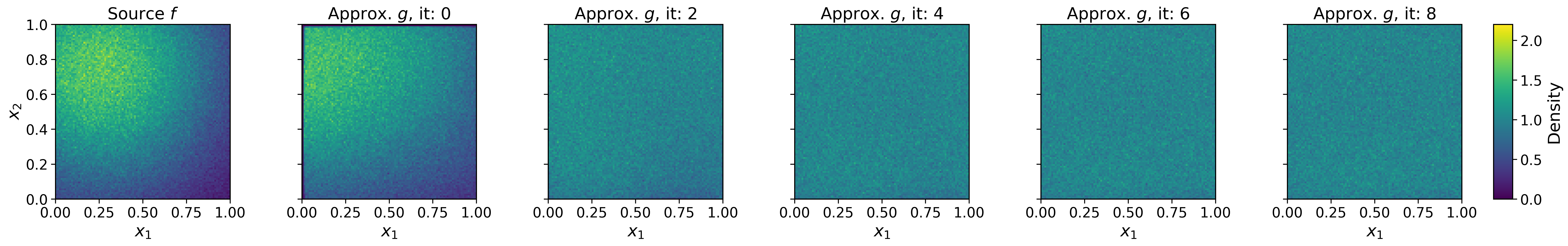}
\caption{Gaussian distribution into uniform.}
\end{subfigure}
\vspace{0.5cm}
\begin{subfigure}{0.95\linewidth}
\centering
\includegraphics[width=0.99\linewidth]{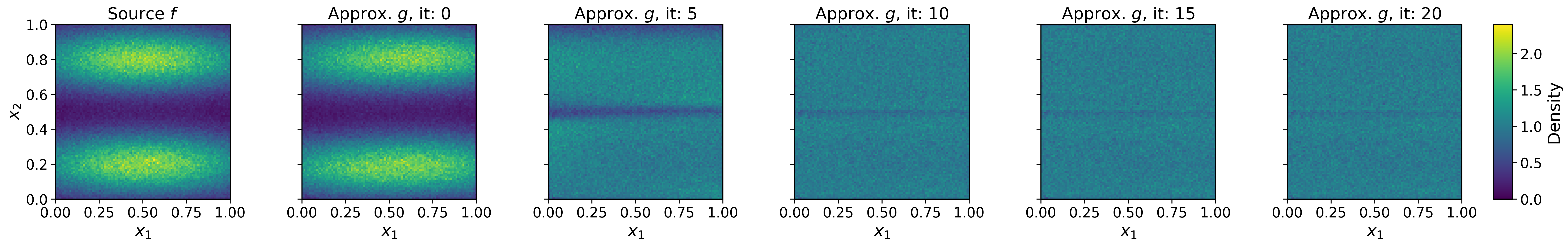}
\caption{Two Gaussian distributions into uniform.}
\end{subfigure}
\vspace{0.5cm}
\begin{subfigure}{0.95\linewidth}
\centering
\includegraphics[width=0.99\linewidth]{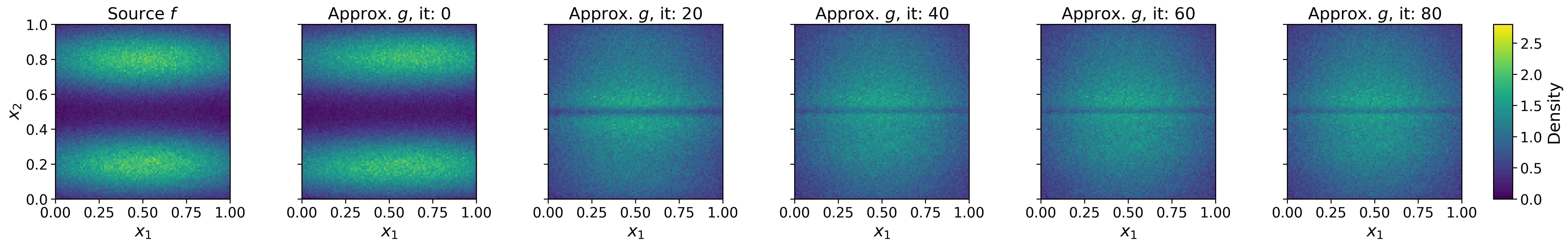}
\caption{Two Gaussian distributions into a Gaussian.}
\end{subfigure}
\caption{
Optimal transport Monge-Amp\`ere problem: transporting gaussian distribution. 
Visualization of the transport of the density function at several iterations of the splitting algorithm with histograms based on $10^6$ sampling points. 
}
\label{fig:ot_gauss_square_comp}
\end{figure}

\Cref{fig:ot_gauss_square_comp} shows histograms obtained from $10^6$ points sampled from the initial densities $f_i$ (left), together with the corresponding transported point distributions at different iterations of the algorithm for the three test cases. In the latter two test cases, we observe a loss of accuracy in the approximation of the transport map along the line $y=0.5$. This region separates the two modes of the bimodal source and appears to be particularly challenging for the numerical approximation.

\section{Conclusions}\label{sec:conc}

We have presented a novel algorithm based on a least-squares approach to solve various second order fully nonlinear equations. 
The approach relies on mathematical programming techniques to treat nonlinearities and a Deep Ritz neural network approach to solve linear variational problems in a more computationally efficient way than traditional finite element techniques. 
The usage of an adaptive sampling strategy in the algorithm allows to significantly accelerates the iterations without sacrificing accuracy. 

Numerical experiments have been proposed along two directions: first we have shown that the method addresses efficiently Dirichlet problems for several 2D fully nonlinear equations. 
Second, we have extended the approach to the optimal transport Monge-Amp\`ere equation, by dealing with transport boundary conditions. The results show that the method can achieve accurate approximations for a range of test cases, while also highlighting limitations for more challenging transport configurations, where localized inaccuracies may persist despite additional splitting iterations. On the other hand, the approach has been shown to be more adapted to strongly nonlinear problems, or solutions with strong gradients, than pure neural networks strategy such as PINNs. 

Future work may include the extension of this approach to higher space dimensions, or the treatment of other nonlinearities.  


\section*{Declaration of generative AI and AI-assisted technologies in the manuscript preparation process.}
During the preparation of this work the authors used ChatGPT in order to help with the plotting and latex editing. After using this tool, the authors reviewed and edited the content as needed and take full responsibility for the content of the published article.

\section*{Declaration of competing interest}
The authors declare that they have no known competing financial interests or personal relationships that could have appeared to influence the work reported in this paper.

\section*{Acknowledgements}
This research did not receive any specific grant from funding agencies in the public, commercial, or not-for-profit sectors.
The authors thank Prof. Marco Picasso (EPFL) for fruitful discussions. 





\bibliographystyle{unsrt}

\bibliography{biblio}
\end{document}